\documentclass{article}

\usepackage{arxiv}

\usepackage{natbib}
\usepackage{doi}

\usepackage{mathtools,amsmath,amssymb,bm,amsthm}

\usepackage{graphicx}
\graphicspath{ {./images/} }
\usepackage[usenames,dvipsnames]{color}
\usepackage{enumitem}
\setlist{nosep}
\usepackage{multirow}
\usepackage{booktabs}
\usepackage{longtable}
\usepackage{lscape,tabularx,makecell}
\usepackage{subfigure}
\usepackage{microtype}
\usepackage{alphalph}
\usepackage{hyperref}

\newcommand{\layer}[1]{\textsf{#1}}

\def\isarxiv{}

\title{Mixed-integer Optimisation of Graph Neural Networks for Computer-Aided Molecular Design}

\author{
Tom McDonald \\
Delft Institute of Applied Mathematics \\
Delft University of Technology \\
Delft, The Netherlands \\
\texttt{thjn.mcdonald@gmail.com}
    \And                
\href{https://orcid.org/0000-0003-2848-2809}{\includegraphics[scale=0.06]{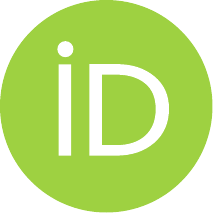}\hspace{1mm}Calvin Tsay} \\
	Department of Computing\\
	Imperial College London\\
	London, UK \\
	\texttt{c.tsay@imperial.ac.uk} \\
	\AND
	\href{https://orcid.org/0000-0001-8885-6847}{\includegraphics[scale=0.06]{orcid.pdf}\hspace{1mm}Artur M. Schweidtmann} \\
Department of Chemical Engineering \\
Delft University of Technology \\
Delft, The Netherlands \\
\texttt{A.Schweidtmann@tudelft.nl} \\
    \And
 \href{https://orcid.org/0000-0002-1814-3515}{\includegraphics[scale=0.06]{orcid.pdf}\hspace{1mm}Neil Yorke-Smith} \\
STAR Lab \\
Delft University of Technology \\
Delft, The Netherlands \\
\texttt{n.yorke-smith@tudelft.nl} \\
}

\date{}

\renewcommand{\shorttitle}{Mixed-integer Optimisation of Trained GNNs}

\hypersetup{
pdftitle={A template for the arxiv style},
pdfsubject={q-bio.NC, q-bio.QM},
pdfauthor={David S.~Hippocampus, Elias D.~Striatum},
pdfkeywords={First keyword, Second keyword, More},
}

\begin{document}
\maketitle

\begin{abstract}
	ReLU neural networks have been modelled as constraints in mixed integer linear programming (MILP), enabling surrogate-based optimisation in various domains
and
efficient solution of machine learning certification problems. However, previous works are mostly limited to
MLPs.
Graph neural networks (GNNs) can learn from non-euclidean data structures such as molecular structures efficiently and are thus highly relevant to computer-aided molecular design (CAMD). We propose a bilinear formulation for ReLU Graph Convolutional Neural Networks and a MILP formulation for ReLU GraphSAGE models. 
These formulations enable solving optimisation problems with trained GNNs embedded to global optimality.
We apply our optimization approach to an illustrative CAMD case study where the formulations of the trained GNNs are used to design molecules with optimal boiling points.
\end{abstract}


\section{Introduction}
\label{sec:intro}

The modelling and designing of molecules have long been an interest to researchers.
The domains where these methods can be applied range anywhere from fuel design, resulting in molecules with decreased emissions, to designing molecules for drug discovery, possibly saving human lives.  Whereas these methods mostly relied on human expertise and experimentation, Computer Aided Molecular Design (CAMD) has become the de facto state of the art \citep{Achenie02:computer}.  CAMD methods for instance pre-screen a large number of molecules, such that the most promising candidates can be investigated for further testing, saving time and resources. 

An early and still established method used for CAMD is the Quantitative Structure Property Relationship (QSPR). With QSPR, chemical descriptors are manually designed and then used to predict chemical properties \cite{de2003qspr,katritzky1995qspr}.  
Many examples exist in which QSPR regressions include constitutional \cite{katritzky1995qspr, ha2005quantitative}, topological \cite{begam2016computer, de2003qspr}, electrostatic \cite{wessel1995prediction, egolf1994prediction}, geometrical \cite{ivanciuc2002quantitative} and quantum-chemical \citep{de2003qspr, hilal2003prediction} descriptors or combinations of these descriptors. 
A drawback of QSPR methods is that they are heavily dependent on the knowledge of researchers to select which chemical descriptors are important. In addition, these methods design a molecule in the descriptor space and thus also give a solution in the descriptor space; mapping this vector back to (existing) molecules is highly problematic. 
Group contribution~\citep{gani1991group,zhang2015generic} presents an alternative family of modelling methodologies, where compounds are represented as a collection of function groups. 
We note that group contribution methods are sometimes classified as a special case of QSPR regression~\citep{alshehri2020deep,gani2019group}.

Often in molecular design, it is interesting to identify a molecule where a desired property is maximised or minimised. There are different methods to achieve this. For instance, mixed integer linear programming (MILP) formulations of QSPR regressions and optimizing them \cite{austin2016computer}. The goal is to optimize an approximation of a particular property and find a corresponding input molecule that corresponds to the maximal value. 
Many efforts also seek to identifying constraints to eliminate infeasible molecules from the optimization search space~\citep{lampe2015computer,zhang2015generic}. 
Owing to the large search spaces and challenging constraints, the resulting optimization problems are often difficult to solve deterministically; moreover, multiple (suboptimal) solutions may be preferred, e.g., to generate a set of promising, ranked candidates.
Therefore, alternative strategies for the optimization step include genetic algorithms~\citep{herring2015evolutionary,venkatasubramanian1994computer,zhou2017hybrid}, or heuristics such as Tabu search~\citep{lin2005computer,mcleese2010design}. 

The advent of machine learning (ML) models and the increased availability of large data sets, resulted in an increased interest in using ML for property prediction tasks~\citep{alshehri2020deep}.  
There have been various applications of machine learning in CAMD.  Most instances use multi-layer perceptrons (MLPs) in the QSPR methods, where linear or polynomial regression methods are replaced by MLPs to perform regression \cite{austin2016computer}.  
More recently, graph neural networks (GNNs) have been developed for non-euclidean input data types such as molecular graphs. 
GNNs are neural networks that learn using a graph as input. Accompanying the spatial graph information, every node in the graph also has an associated feature vector, storing information about that particular node. The information of a node gets passed through an MLP for every node in the network. However, the information that gets passed through the MLP for a node is not only the feature vector of that node but also the feature vectors of the neighbouring nodes in the graph \citep{wu2020comprehensive}. This allows GNNs to take spatial information into consideration when learning non-euclidean data.

The attraction of using GNN in CAMD is that molecules can naturally be represented as graphs. Every atom in a molecule is represented by a node, and the properties of this molecule are stored in the feature vectors associated with the atom-representing nodes. 
Moreover, GNNs preserve invariance of the graph structure, e.g., rotating a molecule does not affect the prediction. 
There have been multiple studies where GNNs have been used to predict properties of molecules (see \citet{wieder2020compact} for an overview). To use these methods in CAMD, just as with the previously-mentioned QSPR methods, one wants to optimize the modelled properties and see which molecule corresponds to this optimized value.  \citet{rittig2022octanefuels} have done exactly that, using Bayesian optimization and a genetic algorithm to optimize the trained GNNs. However, these methods are not deterministic optimization methods. This means that the found solution might be the local maximum of the trained GNN and not the global maximum. In many cases, it is favourable to know with certainty that the found solution is the global optimum. 

GNNs in chemistry can be categorised into three subgroups \citep{wieder2020compact}: (1) Recurrent GNNs (Rec-GNN), (2) Convolutional GNNs (Conv-GNN) and (3) Distinct Graph Neural Network Architectures (Dist-GNN). We will consider the first two subcategories as they are relevant to this paper. 
All of the previously mentioned graph structures have been applied to learning chemical properties. This includes basic Rec-GNNs \cite{lusci2013deep, scarselli2008graph} and gated variants \cite{mansimov2019molecular, withnall2020building, altae2017low, bouritsas2022improving}. Several Conv-GNNs have also found applications in chemistry, such as spectral Conv-GNNs \cite{liao2019lanczosnet,henaff2015deep} and basic \cite{duvenaud2015convolutional, errica2019fair}, attention \cite{hu2019strategies} and general \cite{gilmer2017neural} spatial conv-GNNs.  For a comprehensive overview of molecular property prediction with graph neural networks, see \citet{wieder2020compact}.

Recently, various MILP formulations have been introduced for Rectified Linear Unit (ReLU) MLPs \citep{anderson2020strong,fischetti2018deep,huchette2023deep,tsay2021partition}.  ReLU MLPs are MLPs where each activation function is a piece-wise linear function called the ReLU function. Due to its piece-wise linear nature, the activation function can be expressed with linear programming constraints using big-M constraints. The other functions in a MLP are affine and thus the whole network can be linearised. 
Besides MLP formulations, NNs have also been solved using deterministic global solvers in a reduced space formulation \cite{schweidtmann2019deterministic}. 
The class of MILP problems can be solved to global optimality using commercial solvers. This young research area has been applied to a wide variety of topics like MLP verification \cite{fischetti2018deep, tjeng2017evaluating, Bunel2017, Dutta2017}, compression of MLPs \cite{kumar2019equivalent, serra2020lossless} and using MLPs as surrogate models in linear programming problems \cite{grimstad2019relu, di2022neural, kody2022modeling, yang2021data}. 
We refer the interested reader to \citet{huchette2023deep} for an overview of methodologies and applications. 

The current work, first reported in the master thesis of McDonald \citep{McDonald22:thesis}, is to our knowledge the first MILP formulation of a trained GNN presented in the literature. 
The importance of such a model is that MILP formulations for GNNs can be used in CAMD, where properties of molecules can be modelled using GNNs and then optimized using MILP formulations of these trained GNNs. 
More recent work by \citet{Zhang23:optimizing} develops symmetry-breaking constraints that can reduce the search space for MILP or other optimization strategies. Furthermore, there are broader applications, namely the use of MILP formulations of GNNs for similar applications as MLPs~\citep{huchette2023deep}, e.g., verification of GNNs, lossless compression of GNNs, and using GNNs as surrogate models in optimization problems. 
The latter may be of particular interest for applications such as integrated molecule and process design~\citep{bardow2010continuous}. 

In particular, this current paper considers two GNN architectures.  The first is the Graph Convolutional Neural Network by \citet{kipf2016semi}.  This neural network is one of the earliest GNN and is used often in GNN applications.  The second is the GraphSAGE network by \citet{hamilton2017inductive}, which learns properties of large graph data by sampling the neighbourhood of nodes instead of using information of all neighbouring nodes.

\paragraph{Contributions}
Summarised, this paper adds to the state-of-the-art in the literature as follows:

\begin{itemize}
    \item We propose a mixed integer non-linear programming formulation of the frequently used Graph Convolutional Network model by \citet{kipf2016semi}. 
    \item We propose a mixed integer linear programming formulation of the GraphSAGE model by \citet{hamilton2017inductive}.
    \item We demonstrate the computational performance of our approach on a case study of optimizing the boiling points of molecules modelled with the GraphSAGE and GCN models.
\end{itemize}

\paragraph{Organisation}
Following this introduction, 
Section~\ref{sec:bg} provides technical background, leading to our main contribution of the MI(N)LP formulations of GNNs in Section~\ref{sec:method}.
Section~\ref{sec:results} reports empirical results on a case study.
Section~\ref{sec:discussion} discusses the models and results, and
Section~\ref{sec:conc} concludes.

\ifdefined\isarxiv
\label{Appendix:GCN-full-formulation}
\label{Appendix:extra-results}
\label{Appendix:extra-results-case-study-model-selection}
\label{appendix:extra-runs}
\label{appendix:ga}
\label{appendix:initial-experiments}
\label{appendix:molecule-constraints}
\label{appendix:relu-as-mip}
Appendices are provided in \citet{McDonald22:thesis}.
\else
\fi

\section{Background}
\label{section:background}\label{sec:bg}
This section introduces the terminology for neural networks needed in the remainder of the paper.
We assume the reader has familiarity with mixed integer (linear) programming, referring to \citet{Wolsey20:book} for an introduction.

\subsection{Multilayer Perceptrons}
\label{subsection: ANN-theory}
A feedfoward multilayer perceptron (MLP) consists of consecutive layers of neurons connected through a directed acyclic graph. A neuron in a particular layer receives a weighted signal from the neurons of the previous layer expressed as a real number. Like synapses in the brain, these neurons get activated when the sum of these signals reaches a particular threshold. The result of this system is a neural network that has the ability to emulate complex non-linear relationships. 

In mathematical terms this translates to a neural network $f(x): \mathbb{R}^m \mapsto \mathbb{R}^n$ built of multiple layers $k \in \left\{1, \dots, K \right\}$, including the input layer $k  = 1$, the hidden layers $k = \left\{2, \dots, K-1 \right\}$ and the output layer $k = K$. Each layer contains $n_k$ neurons. Naturally, the input layer has $n_1$ neurons and receives the input vector $x_1 \in \mathbb{R}^{n_1}$ of the function. Every layer $k$ has an associated weight matrix $w^k \in \mathbb{R}^{n_{k} \times n_{k-1}}$ and a bias vector $b^k \in \mathbb{R}^{n_{k}}$ \cite{Goodfellow-et-al-2016}. 

The values associated with neurons in consecutive layers $x_{k} \in \mathbb{R}^{n_{k}}$ are calculated with a propagation function which is a composition of a set of affine functions and non-linear activation functions. This propagation function takes the inputs from real values of the neurons of the previous layer $x_{k-1} \in \mathbb{R}^{n_{k-1}}$. Thus, for the hidden layers $k = \left\{2, \dots, K-1 \right\}$ we have
\begin{equation}
    \label{eq:propagation-function}
    g^k(x^{k-1}) = x^k = \sigma(w^k x^{k-1} + b^k),
\end{equation}
where $\sigma(\cdot)$ is the activation function. Normally, in the last layer $K$ the activation function is absent. 
Completely composed, the neural network $f(x): \mathbb{R}^{n_1} \mapsto \mathbb{R}^{n_K}$ is defined by \citet{Goodfellow-et-al-2016}:
\begin{equation}
    f(x^1) = x^K = (g^K \circ g^{K-1} \circ \dots \circ g^2 \circ g^1)(x_1).
\end{equation}

The activation function, indicated by $\sigma$ in Eq.~\eqref{eq:propagation-function}, is a non-linear function, which allows the neural network to find a non-linear relationship between input and output data. Commonly used activation functions include the sigmoid, tanh, and ReLU functions. The latter will be the main focus for this this contribution.  It is defined as $\sigma(z) = \max\{0,z\}.$

Supervised learning uses paired data, where each data point consists of an input vector $x$, and a desired output $y$. The goal of the learning task is to tune the weights and biases to minimize a loss function, e.g., mean squared error (MSE), for the predictions and target values \cite{Goodfellow-et-al-2016}.

\subsection{MILP Formulations of Multilayer Perceptrons} 
\label{subsection:MILP-ANN}
Exact MILP formulations of NNs
with ReLU activation functions have been proposed. These exact formulations emulate the ReLU operator using binary activation variables and big-M formulations. We refer to \citet{huchette2023deep} for a survey of methods and applications and defer some details to \ref{appendix:relu-as-mip}.

Consider, for each hidden layer $k \in \left\{1, \dots, K - 1\right\}$, the following MLP layer:
\begin{equation}
\label{eq:ANN-layer}
    x^k = \sigma(W^k x^{k-1} + b^k)
\end{equation}
where $\sigma(\cdot) = \operatorname{max}\left\{0,\cdot\right\}$ is the ReLU function, $x^k \in \mathbb{R}^{n_k}$ is the output of layer $k$, $W^k$ and $b^k$ are respectively the found weights and bias of layer $k$. 
This paper considers the linearisation of \eqref{eq:ANN-layer} by \citet{fischetti2018deep}. The output of the affine equations are decoupled in a positive part $x\geq0$ and negative part $s \geq 0$, 
and a binary activation variable $z$ and big-M activation constraints are introduced.
It is assumed that bounds can be found such that $l \leq w^Ty + b \leq u$. For every neuron $j$ layer $k$ of any neural network where the ReLU function is applied the following set of constraints are introduced:
\begin{subequations}
\label{eq:bigM}
\begin{align}
    x_j^k & \leq u_j^k z_j^k \\
    s_j^k & \leq -l_j^k(1-z_j^k) \\
    z_j^k & \in \{0,1\}.
\end{align}
\end{subequations}
The big-M constraints are applied to every node in the network.

\ifdefined\isarxiv
The following states the formulation for a multilayer perceptron with $K$ layers and $n_k$ nodes $j$ per layer. It assumes the final output layer $K$ to be singular and there not to be a ReLU function on that layer. 
\begin{subequations}
\label{eq:ANN-MILP}
\begin{align}
 \quad &\max x_1^K &\\
\text{s.t.}\quad & \bm{W}^K \bm{x}^{K-1} + b^K = x_1^K \\
\label{eq:affine-layerK}
& \bm{W}^k_j \bm{x}^{k-1} + b^k_j = x^k_j - s^k_j & \forall k \in \{1, \dots, K - 1\}, \forall j \in \{1, \dots, n_k\}\\
\label{eq:affine-layersHidden}
& x_j^k \leq u_j^k z_j^k & \forall k \in \{1, \dots, K - 1\}, \forall j \in \{1, \dots, n_k\}\\
& s_j^k \leq -l_j^k(1-z_j^k) & \forall k \in \{1, \dots, K - 1\}, \forall j \in \{1, \dots, n_k\}\\
& x^k_j, s^k_j \geq 0 & \forall k \in \{1, \dots, K - 1\}, \forall j \in \{1, \dots, n_k\}\\
& z^k_j \in \{0,1\} & \forall k \in \{1, \dots, K - 1\}, \forall j \in \{1, \dots, n_k\}\\
& x^0 \in \Omega
\end{align}
\end{subequations}

In this formulation, $\bm{W}^k_j$ is row $j$ of the weight matrix of layer $k$, which naturally has the same dimension as the output $\bm{x}^{k-1}$ of the previous layer. The first input vector is constrained by the input constraints $\Omega$. These are additional input constraints, containing the input bounds, but also other properties which can constrain the input vector, when used in surrogate models for example. 

As noted by \citet{grimstad2019relu}, this is an exact formulation of the ReLU neural network. This means that the above formulation exactly emulates the trained neural network from which the weight matrices $W^k$ and biases $b$ are extracted. For any given input $x_0$ the output of the MILP formulation and the neural net should have the same outcome. The solution which the MILP solver finds also finds consistent solution variables, with the exception of differing $z_j^k$ variables in case the input node $x_j^k$ is $0$.  Note this has no effect on the output, however. 
\else
\fi

\subsection{Graph Neural Networks}
\label{subsection:GNN-theory}

We now turn from MLPs and `regular' neural networks to GNNs.

\subsubsection{General Graph Neural Network Architecture}
When using GNNs for property prediction, each data point consists of the structure of a graph $G = (V,E)$ represented by the adjacency matrix $A \in \mathbb{R}^{N \times N}$, and properties of the graphs. The properties of these graphs are stored in node feature vectors $X \in \mathbb{R}^{N \times F}$, and can sometimes include edge feature vectors. For our purposes of CAMD, node features will suffice. Every node $i \in V$ has an accompanying feature vector $X_i \in \mathbb{R}^F$. These feature vectors store information about the node in question. In a supervised setting, the data is thus of the form $((X,A) ,y)$. 

Graph convolutional neural networks are divided in spectral and spatial based methods. Spectral based methods are graph neural networks based on graph signal filters. Spatial based methods are generally GNNs consisting of a function which aggregates neighbourhood information and some sort of propagation function, similar to those found in MLPs. The aggregation function sums the the feature vectors of  neighbouring nodes of a node $i$, which is used as input of an affine function. Thereafter, the affine combination of the aggregated feature vectors is passed through an activation function, similar to the feedfoward neural network architecture.  Doing this for every node in the graph constitutes one convolutional layer. After one convolutional layer, every node has a new feature vector. 

Stacking multiple convolutional layers consecutively allows a node $i$ to not only process node feature vector information of its neighbouring nodes $\mathcal{N}(i)$, but also of the neighbours $\mathcal{N}(s)$ of these neighbours $\forall s \in \mathcal{N}(i)$. This works as follows: in the first convolutional layer, for every node $i$, all neighbourhood information is aggregated. In the next layer this is repeated; however, all neighbours of node $i$ have already processed the information of their respective neighbours. This means $i$ also internalises the information of all neighbours removed with a $2$-length path. After $k$ convolutions, node $i$ processes information from all nodes $k$-length paths removed. 

In the following subsections we will discuss the graph aggregation functions of two GNNs, for they define the architectures of the GNNs we consider.

\begin{figure}[tb]
    \centering
    \includegraphics[width=0.8\textwidth]{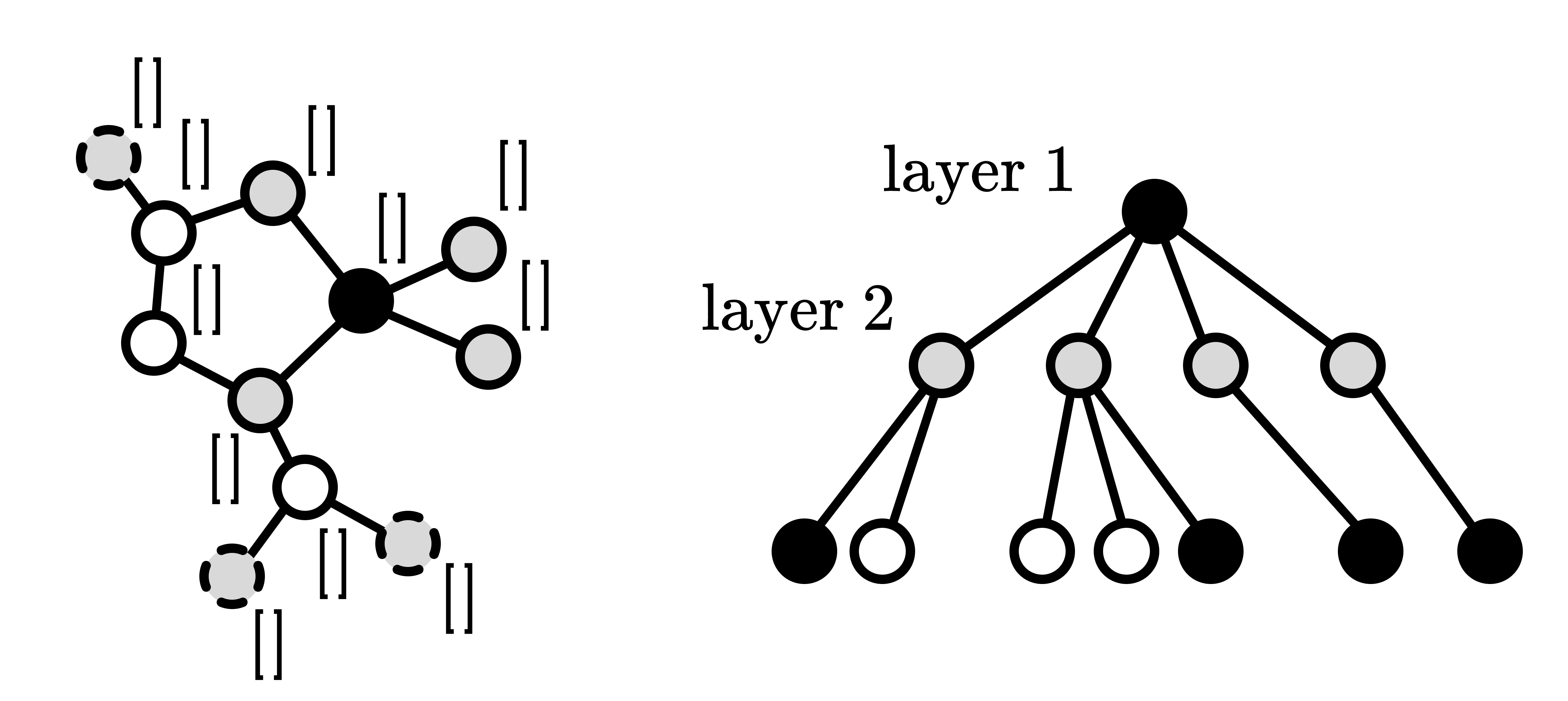}
    \caption{A graph (left) with a feature vector on every node. Neighbourhoods of a node with multiple convolutional layers in a spatial GNN (right)}
    \label{fig:my_label}
\end{figure}

\subsubsection{Graph Convolutional Neural Network (GCN)}
We focus on spatial Conv-GNN methods, which are conceptually similar to non-graph based convolutional neural nets (CNN), as ``spatial-based graph convolutions convolve the central node’s representation with its neighbors’ representations to derive the updated representation for the central node'' \cite{wu2020comprehensive}. \citet{Micheli2009} introduced these spatial graph Neural Networks. Thereafter, many varieties of spatial Graph Neural networks have been introduced. Basic models include PATCHY-SAN, LGCN and GraphSAGE \cite{niepert2016learning, gao2018large, hamilton2017inductive}. All use a combination of convolutional operators, combined with different neighbour selection systems and different aggregators. There is also a set of attention-based spatial approaches which assign different weights for different neighbours to minimise noise \cite{velivckovic2017graph, zhang2018gaan}. Finally, there are more general frameworks which try to unify multiple models in a single formulation as an abstraction over multiple GNNs \cite{monti2017geometric, gilmer2017neural, battaglia2018relational}.

The first GNN we consider is the Graph Convolutional Neural Network (GCN) \cite{kipf2016semi}. It is one of the earlier models which can be considered as a spatial GNN method. The GCN has its roots in spectral graph GNNs as it is a first order Chebychev approximation of the ChebNet \cite{defferrard2016convolutional} architecture, which is a spectral based method. However, this first order approximation is basically a spatial based method. 

\citet{kipf2016semi} introduce the $k$-th convolutional layer in the GCN can be expressed as follows:
\begin{equation}
    H^{(k+1)} = \sigma(\Tilde{D}^{-\frac{1}{2}}\Tilde{A}\Tilde{D}^{-\frac{1}{2}} H^{(k)} W^{(k)}) \label{eq:ConvGNN}
\end{equation}
Here, $\tilde{A}=A+I_{N}$ is the adjacency matrix of the undirected graph $\mathcal{G}$ with added self-connections. $I_{N}$ is the identity matrix, $\tilde{D}_{i i}=\sum_{j} \tilde{A}_{i j}$ and $W^{(k)}$ is a layer-specific trainable weight matrix. $\sigma(\cdot)$ denotes an activation function, such as the $\operatorname{ReLU}(\cdot)=\max (0, \cdot)$. $H^{(k)} \in \mathbb{R}^{N \times n_k}$ is the matrix of activations in the $k^{\text {th }}$ layer; $H^{(0)}=X$, where $X$ is a matrix of node feature vectors $X_i$ belonging to node $i$ in the graph. 

The formula to find feature $j$ for a node $i$ in layer $k + 1$ shows the spatial nature of the GCN network:
\begin{equation}
    H_{ij}^{(k+1)} = \sigma \left( \left ( W^{(k)}_j \right )^T \sum_{l \in N^{+}(i)} \frac{1}{\sqrt{d^+(i)}\sqrt{d^+(l)}}  \left ( H^{(k)}_l \right )\right )
\end{equation}
Here, $\mathcal{N}^+(i)$ is the neighbourhood of $i$ including $i$ itself, and $d^+(i)$ is the degree of node $i$. The aggregation function is a normalised sum of all the feature vectors of $l \in \mathcal{N}^+(i)$ in layer $k$. Thereafter, just as in MLP models, an affine combination is taken of the aggregated feature vectors and passed through an activation function $\sigma$.

\subsubsection{GraphSAGE Network}

The GraphSAGE network is another spatial convolutional neural network, developed to learn large graph networks. Its input is merely one large graph $G = (V,E)$ on which it performs the learning task. When trained, the GraphSAGE network can classify nodes, without having seen all nodes of the network. This means that it can generalise to unseen nodes in the network. 

GraphSAGE also uses an aggregation scheme, for instance the $\verb|mean|$, $\verb|max|$, $\verb|ltsm|$ or $\verb|add|$ aggregation scheme. However, for node $i$, GraphSAGE does not aggregate over all feature vector of its neighbours $\mathcal{N}(i)$, but over a randomised subset of the neighbourhood. This allows it to learn large graphs. The aggregated subset of the neighbourhood vectors gets concatenated with the vector of the root node $i$. The concatenated vectors are then multiplied with a learned weight matrix $W^k \in \mathbb{R}^{(n_{k+1} \times 2n_k )}$ which consecutively passes through an activation function $\sigma$. Finally, the vector gets normalised. As usual, all previously described steps are performed for all nodes $i \in V$. 

For this paper we are interested in the GraphSAGE network as it is linearisable, when specific choices for the hyper-parameters are made. We set the sampling to select all neighbours with a probability of $1$.  This can be interpreted such that we don't have a sampling function. The chosen activation function is the ReLU function. We choose the aggregate scheme to be $\verb|add|$, which means that we add all feature vectors of the neighbouring nodes. The propagation function becomes:
\begin{equation}
\label{eq:GraphSAGE-layer-original}
H^{(k+1)}_i=\sigma\left(H^{(k)}_i \cdot W_{1}^{(k)}+\sum_{l \in \mathcal{N}(i)} H^{(k)}_l \cdot W_{2}^{(k)}\right)
\end{equation}

The matrices $W_{1}^{(k)}, W_{2}^{(k)} \in \mathbb{R}^{(n_{k+1} \times n_k)}$ are a split representation of the matrix $W^k \in \mathbb{R}^{(n_{k+1} \times 2n_k )}$, introduced for legibility. Using the add function is also more natural when predicting the boiling points for chemical compounds, which we will discuss in the next subsection.

\section{Methods}
\label{section:methods}\label{sec:method}

This section provides the main contribution of the paper, by presenting the novel formulations of graph neural networks as MI(N)LPs, starting from the multilayer perceptron MILP formulation of Section~\ref{subsection:MILP-ANN}.  Section~\ref{subsection:linear-GNN} describes a formulation of the Graph Convolutional Network, which is linear (MILP) for a fixed graph structure and bi-linear (MINLP) for a variable graph structure. Then, Section~\ref{subsection:GraphSAGE} describes a MILP formulation of the GraphSAGE architecture.
Section~\ref{subsection:input-constraints-simple} shows how to add domain-specific background knowledge (inductive bias).
Finally, Section~\ref{sec:method:fbbt} applies bound tightening techniques to the formulations.

\subsection{MI(N)LP Formulation for GCN Models}
\label{subsection:linear-GNN}

We wish to train a graph neural network that finds the function $f: \{0,1\}^{|N| \times |N|} \times \mathbb{R}^{|N| \times |F|} \mapsto \mathbb{R}$.  This function maps an adjacency matrix $A$ and a feature vector $X$, with $|F|$ features, to a singular output. The function is a composition of GNN layers, a pooling layer, and MLP layers, where the latter takes a fingerprint of the graph and maps it to a singular output.  Mathematically this constitutes:
\begin{equation}
\label{eq:neural-net-conceptual}
    f(A, X) = \layer{MLP }(\layer{POOL }(\layer{GNN}(A,X)))
\end{equation}
where \layer{POOL} is the pooling layer.  This subsection states the linear MILP formulations for graph neural networks in case the graph structure is predetermined.

\subsubsection{Predetermined Graph Structure}
\label{subsubsection:lin-known-GNN}
In the following subsection the Graph Convolutional Network layers as described in Section~\ref{subsection:GNN-theory} are formulated as an MILP.  
We first consider the GCN, which has the layer-wise propagation rule defined by Eq.~\eqref{eq:ConvGNN}. 
To linearise Eq.~\eqref{eq:ConvGNN} we employ the big-M formulation stated in Eq.s~\eqref{eq:ANN-MILP}.  As described in Section~\ref{subsection:MILP-ANN}, the ReLU constraints are the same for every node, with altering big-M values (lower and upper bounds).  For the MLP structure, there were $K$ layers with $n_k$ neurons in each $k^{\text{th}}$ layer. For the GCN structure we have the same but for every node $i \in \{1, \dots, N\}$.  In the first layer these input nodes are the feature vectors for those nodes. 

Since the ReLU constraints stay the same, merely the left hand side of Eq.~\eqref{eq:affine-layerK} and Eq.~\eqref{eq:affine-layersHidden} need to be altered to represent the linear part of the GCN layer as described between the brackets in Eq.~\eqref{eq:ConvGNN}. To simplify the formulation we write $\Bar{A} = \tilde{D}^{-\frac{1}{2}} \tilde{A} \tilde{D}^{-\frac{1}{2}}$. The entries of this matrix are the following:
\begin{equation}
\label{eq:relu-logic-norm-term}
    \Bar{A}_{il} =
    \begin{cases}
      0, & \text{if $\tilde{A}_{il} = 0$}\\
      \frac{1}{\sqrt{d^+(i)}\sqrt{d^+(l)}}, & \text{if $\tilde{A}_{il} = 1$}
    \end{cases}       
\end{equation}
where $d^+(i)$ is the cardinality of the adjacent set $N^+(i)$ for node $i$, where the $+$ indicates that it also includes self loops. 

\citet{kipf2016semi} note the following, with $N$ nodes in our graph
\begin{equation}
    Y^0 = X = 
    \begin{bmatrix}
    X_1 \\
    \vdots \\
    X_N \\
    \end{bmatrix}
\end{equation}
where $X_i$ is a row vector containing the features of node $i$. $H_{ij}^{(k)}$ is considered, which is the $j$-th neuron of node $i$, after $k$ GCN layers. The value of this neuron is found as follows (the activation function $\sigma$ is omitted from every line for clarity):
\begin{subequations}
{\allowdisplaybreaks
\begin{align}
    H^{(k)}_{ij} &= \left ( \Bar{A} H^{(k-1)} W^{(k)} \right )_{ij} \\
    &= \Bar{A_i} \begin{bmatrix}
        H^{(k-1)}_1 W^{(k)}_j \\
        \dots \\
        H^{(k-1)}_N W^{(k)}_j \\
    \end{bmatrix} \\
    &= \sum_{l \in N^{+}(i)} \frac{1}{\sqrt{d^+(i)}\sqrt{d^+(l)}} \left ( W^{(k)}_j \right )^T \left ( H^{(k-1)}_l \right )^T
\end{align}}
\end{subequations}

It is commonplace to write vectors in column notation for linear programming, so we will deviate from \citet{kipf2016semi}, replacing $\left(H_{l}^{(k-1)}\right)^T$ by $\left(H_{l}^{(k-1)}\right)$.  The MILP formulation becomes: 
\begin{subequations}
\label{eq:GCN-MILP}
\begin{multline}
    \sum_{l|\tilde{A}_{il} = 1} \frac{1}{\sqrt{d^{+}_i d^{+}_l}} {\bm{W}^{k}_j}^T \bm{H}^{(k-1)}_l = H_{ij}^k - S_{ij}^k \\ \forall k \in \{1, \dots, K\}, \forall i \in \{1, \dots, N\}, \forall j \in \{1, \dots, n_k\}    \label{eq:conditional-hidden-layer}
\end{multline}
{\allowdisplaybreaks
    \begin{align}
    & H_{ij}^k \leq U_{ij}^k Z_{ij}^k & \forall k \in \{1, \dots, K\}, \forall i \in \{1, \dots, N\}, \forall j \in \{1, \dots, n_k\}\\
    & S_{ij}^k \leq -L_{ij}^k(1-Z_{ij}^k) & \forall k \in \{1, \dots, K\}, \forall i \in \{1, \dots, N\}, \forall j \in \{1, \dots, n_k\} \label{eq:GCN-MILP-c}\\
    & 0 \leq H_{ij}^k, S_{ij}^k \in \mathbb{R} & \forall k \in \{1, \dots, K\}, \forall i \in \{1, \dots, N\}, \forall j \in \{1, \dots, n_k\}\\
    & Z_{ij}^k \in \{0,1\} & \forall k \in \{1, \dots, K\}, \forall i \in \{1, \dots, N\}, \forall j \in \{1, \dots, n_k\}\\
    & A,H^0 \in \Omega
    \end{align}}
\end{subequations}
where $i$ indicates node $i$, $j$ feature $j$, and $k$ the corresponding GCN layer. $d^{+}_i$ represents the degree $+ 1$ of node $i$.

In the above, $H^0 \in \Omega$ indicates a restriction of the input space (see Section~\ref{subsection:input-constraints-simple}).  For a molecule we can see these as constraints as providing domain knowledge such that only physically possible molecules are considered in the input space. One example could be that node $i$, has a feature $x_{ij} = 1$ indicating that it is a carbon molecule. In that case $\sum_{j} A_{ij} < 5$, meaning it cannot share a bond with more than $4$ other molecules since the amount of atomic bonds is maximally 4 for a carbon molecule. 

Combining all constraints of \eqref{eq:GCN-MILP} with \eqref{eq:ANN-MILP} almost finalises our formulation of function \eqref{eq:neural-net-conceptual} in case the graph structure $A$ is known. We still need to include a pooling layer. The pooling layer is a function which operates over all neuron outputs of the nodes of the graph to combine them into a single vector. There are multiple options for pooling layers but we utilise a sum pooling layer to maintain a linear model: 
\begin{equation}
\label{eq:pooling-layer}
    x^0_j = \sum_{i = 1}^{N} H_{ij}^K \quad \forall j \in \{1, \dots, n_K\}
\end{equation}

Naturally, the amount of entries of the input vector $x_0$ of the MLP layer, must match the number of neurons in the final layer $K$ of the GCN layers. 

\subsubsection{Variable Graph Structure}
\label{subsubsection:linearising-non-linear-terms}

The formulation of the previous subsection is linear in case the structure of the graph is known. If $A$ is unknown, this formulation is non-linear. There are many examples where we wish to find the graph structure accompanying an optimal solution. The non-linear terms in the above formulation, in case the structure is unknown, are $\sum_{l|\tilde{A}_{il} = 1}$ and $\frac{1}{\sqrt{d_i^+d_l^+}}$, where $l$ in the latter is dependant on the former non-linear term. 
We describe a bi-linear formulation for these terms, which MIQP/MIQCP solvers such as Gurobi can accomodate.

\paragraph{A linear conditional sum}
The sum in the formulation is conditional and thus non-linear. To linearise the sum a support variable $\bm{b}^k_{il}$ is introduced. This new support variable follows the following logic for two nodes $i$ and $j$:
\begin{equation}
\label{eq:sum-support-logic}
    \bm{b}^k_{il} =
    \begin{cases}
      0 & \text{if $\tilde{A}_{il} = 0$}\\
      \bm{H}^k_l & \text{if $\tilde{A}_{il} = 1$}
    \end{cases}       
\end{equation}
If this logic is implemented the sum over all $\bm{b}^k_{il}$ results in the same outcome as the conditional sum. For every node $i$, $\bm{b}^k_{il}$ are only equal to the output of ReLU layer if they are connected to node $i$. 

The logic as described in \eqref{eq:sum-support-logic} can be implemented 
\ifdefined\isarxiv
\else
in the same way as was done in Eq.~\eqref{eq:relu-logic} 
\fi
by using big-M constraints, because the entries of $\tilde{A}$ are binary.  Replacing Eq.~\eqref{eq:conditional-hidden-layer} by the following constraints for $ k \in \{1, \dots, K\}, i \in \{1, \dots, N\}, j \in \{1, \dots, n_k\}$ we have removed the conditional sum: 
\begin{subequations}
\label{eq:sqrt-hidden-layer}
\begin{align}
    \label{eq:sqrt-hidden-layer-a}
    \sum_{l} \frac{1}{\sqrt{d^{+}_i d^{+}_l}} {\bm{W}^{k}_j}^T \bm{b}^{(k-1)}_{il} &= H_{ij}^k - S_{ij}^k \\
    \label{eq:sqrt-hidden-layer-b}
    \bm{H}^{(k-1)}_l - \bm{M}(1 - \tilde{A}_{il}) \leq \bm{b}^{(k-1)}_{il} &\leq \bm{H}^{(k-1)}_l + \bm{M}(1 - \tilde{A}_{il})& \\
    \label{eq:sqrt-hidden-layer-c}
    -\bm{M}(\tilde{A}_{il}) \leq \bm{b}^{(k-1)}_{il} &\leq \bm{M}(\tilde{A}_{il}) 
\end{align}
\end{subequations}
In case node $i$ is not connected to node $l$, then $\tilde{A}_{il} = 0$. In that case constraint \eqref{eq:sqrt-hidden-layer-c} forces $\bm{b}^{(k-1)}_{il} = 0$. In case both nodes are connected, $\tilde{A}_{il} = 1$ and $\bm{b}^{(k-1)}_{il}$ is constrained by \eqref{eq:sqrt-hidden-layer-b} such that it is equal to $\bm{H}^{(k-1)}_l$. 

\paragraph{A linear normalisation term} 
We are still left with $\frac{1}{\sqrt{d_i^+d_l^+}}$, which is also non-linear. The term is also multiplied with the variable vector $b_{il}^{k}$, which makes the entire constraint non-linear. It is possible to remove the fraction and the square root, albeit in a contrived way, adding a lot of extra variables. We note that the cardinality of the co-domain of the function $g(i,l) = \frac{1}{\sqrt{d_i^+d_l^+}}$, is upper bound by the maximum degree of the graph $d_{max}$, adding $1$ for the self loops. In case of molecules this is $4$ for instance. This means that the function $g$ has a maximum of $(4+1)^2$ outcomes. We can index these outcomes in a $(d_{max} + 1)^2$ long vector $g$, where at index $p = d^+_i (d_{max} + 1) + d^+_l$, $g_p = \frac{1}{\sqrt{d_i^+d_l^+}}$. The function $g(i,l)$ is undefined in case $d_i^+ = 0$ or $d_l^+ = 0$. In these cases $g_p = 0$.

Using linear constraints, we can linearise the fractional term in Eq.~\eqref{eq:sqrt-hidden-layer} by the following set of equations
\begin{subequations}
\label{eq:normalisation-term-linearisation}
\begin{multline}
    \label{eq:normalisation-term-linearisation-a}
    \sum_{l} s_{il} {\bm{W}^{k}_j}^T \bm{b}^{(k-1)}_{il} = H_{ij}^k - S_{ij}^k \\ \forall k \in \{1, \dots, K\},\forall i \in \{1, \dots, N\}, \forall j \in \{1, \dots, n_k\}
\end{multline}
{\allowdisplaybreaks
\begin{align}
    &d_i^+ = \sum_{j} \tilde{A}_{ij} &\forall i \in \{1, \dots, N\} \\
    &p_{il} = d^+_i (d_{max} + 1) + d^+_l & \forall i,l \in \{1, \dots, N\}\\
    &0 = p_{il} - 1c^{il}_1 - 2c^{il}_2 - \dots \\
        &\qquad - (d_{max} + 1)^2 c^{il}_{(d_{max} + 1)^2} & \forall i,l \in \{1, \dots, N\}\\
    &1 = c^{il}_1 + \dots + c^{il}_{(d_{max} + 1)^2}& \forall i,l \in \{1, \dots, N\} \\
    &c^{il} \in \{0,1\}^{(d_{max} + 1)^2} & \forall i,l \in \{1, \dots, N\}\\
   & s_{il} = c^{il}_1 g_1 + \dots + c^{il}_{(d_{max} + 1)^2} g_{(d_{max} + 1)^2}& \forall i,l \in \{1, \dots, N\}
\end{align}}
\end{subequations}

With this set of equations, we are mapping the index $p_{il}$ to its corresponding value in vector $g$, which is a set of predetermined parameters. Since the structure of graph stays the same over all GCN layers, we only have to add these constraints once and not for every layer $k$. 
The resulting model is bi-linear as it involves multiplication of decision variables $s_{il}$ and $b_{il}$.

\subsubsection{Full MINLP Formulation of the GCN GNN}
\label{subsection:non-linear-GNN}

\label{subsubsection:single term bi linear}
Section~\ref{subsubsection:linearising-non-linear-terms} presents a reformulation for the non-linear terms of Eq.~\eqref{eq:conditional-hidden-layer}.
Specifically, the conditional sum is reformulated as Eq. \eqref{eq:sqrt-hidden-layer} and the normalisation term as Eq. \eqref{eq:normalisation-term-linearisation}, resulting in an overall bi-linear formulation.  
For GCN models, we can combine the binary variables introduced by Eq. \eqref{eq:sqrt-hidden-layer} and Eq. \eqref{eq:normalisation-term-linearisation}
Notice how for every layer $k$, Eqs. \eqref{eq:sqrt-hidden-layer-b} and \eqref{eq:sqrt-hidden-layer-c} constrain whether or not the feature vector of the neighbours of node $i$ are included in the conditional sum. We can simplify this by incorporating it in the variable which encompasses the linearised normalisation term $s_{il}$. Once again we incorporate the following logic with big-M constraints for $i,l \in \{1, \dots, N\}$ :
\begin{equation}
\label{eq:bi-linear-simplification logic}
    \hat{s}_{il} =
    \begin{cases}
      0 & \text{if $\tilde{A}_{il} = 0$}\\
      s_{il} & \text{if $\tilde{A}_{il} = 1$}
    \end{cases}       
\end{equation}

This enforces that the feature vector of a neighbouring node of $i$ is only included if $\tilde{A}_{il} = 1$, which is the same as $\sum_{l|\tilde{A}_{il} = 1}$. The resulting bi-linear MINLP formulation becomes:
\begin{subequations}
\label{eq:single-bilinear-constraints}
\begin{multline}
    \label{eq:single-bilinear-constraints-a}
    \sum_{l} \hat{s}_{il} {\bm{W}^{k}_j}^T \bm{H}^{(k-1)}_l = H_{ij}^k - S_{ij}^k \\ \forall k \in \{1, \dots, K\}, \forall i \in \{1, \dots, N\}, \forall j \in \{1, \dots, n_k\} 
\end{multline}
{\allowdisplaybreaks
\begin{align}
    \label{eq:single-bilinear-constraints-b}
    &d_i^+ = \sum_{j} \tilde{A}_{ij} & \forall i \in \{1, \dots, N\}\\
    &p_{il} = d^+_i (d_{max} + 1) + d^+_l & \forall i,l \in \{1, \dots, N\}\\
    &0 = p_{il} - 1c^{il}_1 - 2c^{il}_2 - \dots \\
        &\qquad - (d_{max} + 1)^2 c^{il}_{(d_{max} + 1)^2} & \forall i,l \in \{1, \dots, N\}\\
    &1 = c^{il}_1 + \dots + c^{il}_{(d_{max} + 1)^2} & \forall i,l \in \{1, \dots, N\}\\
    &c^{il} \in \{0,1\}^{(d_{max} + 1)^2} & \forall i,l \in \{1, \dots, N\}\\
    \label{eq:single-bilinear-constraints-g}
    &s_{il} = c^{il}_1 g_1 + \dots + c^{il}_{(d_{max} + 1)^2} g_{(d_{max} + 1)^2}& \forall i,l \in \{1, \dots, N\}\\
    \label{eq:single-bilinear-constraints-h}
    &s_{il} -M (1-A_{il}) \leq \hat{s}_{il} \leq M (1- A_{il}) + s_{il}& \forall i,l \in \{1, \dots, N\}\\
    \label{eq:single-bilinear-constraints-i}
    &-M A_{il} \leq \hat{s}_{il} \leq M A_{il}& \forall i,l \in \{1, \dots, N\}
\end{align}}
\end{subequations}

In this formulation, constraints \eqref{eq:single-bilinear-constraints-b}--\eqref{eq:single-bilinear-constraints-g} describe the linearisation of the normalisation term as described in Section~\ref{subsubsection:linearising-non-linear-terms} and constraints \eqref{eq:single-bilinear-constraints-h} and~\eqref{eq:single-bilinear-constraints-i} incorporate the conditional-sum logic from Eq.~\eqref{eq:bi-linear-simplification logic}.  Note that we only have to find $\hat{s}_{il}$ once for every layer, since the structure of the molecule is constant per layer. 

A reader might note that when the degree of either node $i$ or node $j$ is zero, this means that $s_{il}$ will automatically be equal to zero, and thus the introduction of Eqs. \eqref{eq:single-bilinear-constraints-h} and \eqref{eq:single-bilinear-constraints-i} might be superfluous. However, there could be an instance when both $i$ and $l$ have a degree higher than $0$, but still not be connected. In that case $s_{il}$ is not zero, and thus the extra constraints need to be introduced.

\subsection{GraphSAGE}
\label{subsection:GraphSAGE}

So far we have successfully formulated an exact MINLP (bi-linear) representation of a GCN.
This section describes an MILP formulation for a more recent and popular GNN architecture, the \textit{GraphSAGE} model by \citet{hamilton2017inductive}.

In this paper the activation function and affine layer are described by the following equation:
\begin{equation}
\label{eq:GraphSAGE-layer-original2}
f^{(t)}(v)=\sigma\left(f^{(t-1)}(v) \cdot W_{1}^{(t)}+\sum_{w \in N(v)} f^{(t-1)}(w) \cdot W_{2}^{(t)}\right)
\end{equation}
where $f^{(t)}(v)$ describes the feature vector of node $v$ after $t$ GraphSAGE layers and $\sigma$ describes an activation function, which for this paper will once again be the ReLU activation function. 

As can be seen from Eq.~\eqref{eq:GraphSAGE-layer-original2}, after training a neural net with $K$ GraphSAGE layers, it finds two weight matrices for every layer $t$. The first weight matrix $W_{1}^{(t)}$, which we will refer to as the root weight, is multiplied with the feature vector of the previous layer $f^{(t-1)}(v)$. The second weight matrix $W_{2}^{(t)}$ is multiplied with the neighbouring feature vectors of node $v$. In case the adjacency matrix of the graph is unknown, this neighbourhood of $v$ is a non-linear relation. The rest of the model is linear.

For consistency, we rewrite Eq.~\eqref{eq:GraphSAGE-layer-original2} in a notation similar to the presentation in Section \ref{subsection:linear-GNN}. We find the following for node $i \in \left\{1, \dots, N \right\}$, feature $j \in \left\{1, \dots, n_k \right\}$ and layer $k \in \left\{1, \dots, K \right\}$:
\begin{equation}
    \label{eq:GraphSAGE-layer-rewritten}
     {H}^{(k)}_{ij} =
     \sigma \left( 
     { {\bm{{\hat{W}}}}^{k}_{j} }^{T} \bm{H}^{(k-1)}_{i} + 
     { {\bm{{\bar{W}}}}^{k}_j }^{T} \sum_{l | A_{il} = 1 } 
     \bm{H}^{(k-1)}_{l} 
     \right) 
\end{equation}
where $H^{(k)}_{ij} \in \mathbb{R}$ is the feature $j$ of node $i$ after $k$ layers, and $A_{il}$ is the adjacency matrix without self loops.

To remove the conditional sum we again introduce big-M constraints and support variables to encode the following logic:
\begin{equation}
\label{eq:sum-support-logic-GraphSAGE}
    \bm{b}^k_{il} =
    \begin{cases}
      0 & \text{if ${A}_{il} = 0$}\\
      \bm{H}^k_l & \text{if ${A}_{il} = 1$}
    \end{cases}       
\end{equation}

The full MILP formulation including the pooling layer becomes:
\begin{subequations}
\allowdisplaybreaks
\label{eq:GraphSAGE-MILP}
\begin{multline}
    \label{eq:graph-conv-relu-1}
    ({\bm{\hat{W}}^{k}_j})^T \bm{H}^{(k-1)}_{i} + ({\bm{\bar{W}}^{k}_j})^T \sum_{l} \bm{b}^{(k-1)}_{il} = H_{ij}^k - S_{ij}^k \\ \forall k \in \{1, \dots, K\}, \forall i \in \{1, \dots, N\}, \forall j \in \{1, \dots, n_k\}\\
\end{multline}
{\allowdisplaybreaks
    \begin{align}
        \label{eq:graph-conv-relu-2}
        & H_{ij}^k \leq U_{ij}^k Z_{ij}^k & \forall k \in \{1, \dots, K\}, \forall i \in \{1, \dots, N\}, \forall j \in \{1, \dots, n_k\}\\
        \label{eq:graph-conv-relu-3}
        & S_{ij}^k \leq -L_{ij}^k(1-Z_{ij}^k) & \forall k \in \{1, \dots, K\}, \forall i \in \{1, \dots, N\}, \forall j \in \{1, \dots, n_k\}\\[2ex]
        \label{eq:graph-conv-conditional-support-1}
        &\bm{H}^{k}_l - \bm{M}(1 - {A}_{il}) \leq \bm{b}^{k}_{il} & \forall k \in \{0, \dots, K - 1\}, \forall i ,l \in \{1, \dots, N\} \\
        \label{eq:graph-conv-conditional-support-2}
        & \bm{b}^{k}_{il} \leq \bm{H}^{k}_l + \bm{M}(1 - {A}_{il}) & \forall k \in \{0, \dots, K - 1\}, \forall i ,l \in \{1, \dots, N\} \\
        \label{eq:graph-conv-conditional-support-3}
        &-\bm{M}({A}_{il}) \leq \bm{b}^{k}_{il} \leq \bm{M}({A}_{il}) & \forall k \in \{0, \dots, K - 1\}, \forall i ,l \in \{1, \dots, N\}  \\
        \label{eq:graph-conv-conditional-support-4}
        & \bm{H}^{0}_i = \bm{x}_i & \forall i \in \{1, \dots, N \} \\[2ex]
        \label{eq:graph-conv-pooling}
        & {H}^{*K}_j = \sum_i H_{ij}^K & \forall j \in \{1, \dots, n_K \} \\[2ex]
        & 0 \leq H_{ij}^k, S_{ij}^k \in \mathbb{R} & \forall k \in \{1, \dots, K\}, \forall i \in \{1, \dots, N\}, \forall j \in \{1, \dots, n_k\}\\
        & Z_{ij}^k \in \{0,1\} & \forall k \in \{1, \dots, K\}, \forall i \in \{1, \dots, N\}, \forall j \in \{1, \dots, n_k\}\\
        \label{eq:GraphSAGE-MILP-formulation-input-constraints}
        & \bm{x}, A \in \Omega
    \end{align}}
\end{subequations}

Intuitively, Eqs.~\eqref{eq:graph-conv-relu-1}--\eqref{eq:graph-conv-relu-3} reformulate the ReLU function in Eq.~\eqref{eq:GraphSAGE-layer-rewritten}, Eqs.~(\ref{eq:graph-conv-conditional-support-1}--\ref{eq:graph-conv-conditional-support-3}) are the big-M constraints to enforce the logic in Eq.~\eqref{eq:sum-support-logic-GraphSAGE}. The values of these big-M constraints are the same as the upper bounds $U_{ij}^k$ (as explained in Section~\ref{subsubsubsection:FBBT-GraphSAGE}).  Eq.~\eqref{eq:graph-conv-conditional-support-4} defines the input feature vector $\bm{x}_i$ of node $i$ and Eq.~\eqref{eq:graph-conv-pooling} is the sum pooling layer in layer $K$.  Finally, Eq.~\eqref{eq:GraphSAGE-MILP-formulation-input-constraints} represent the input constraints as described in Section~\ref{subsection:input-constraints-simple}.  

Note that a drawback of this method compared to the MINLP formulation is that constraints \eqref{eq:graph-conv-conditional-support-1}, \eqref{eq:graph-conv-conditional-support-2} and~\eqref{eq:graph-conv-conditional-support-3} are calculated for every layer $k$. This increases the number of constraints significantly (by $\mathcal{O}(n^2 n_k)$ constraints per layer $k$). For the GCN network this is only $\mathcal{O}(n^2)$. For networks where there are a lot of nodes per hidden layer, the GraphSAGE network will have significantly more constraints than the GCN network.

\subsection{Constraining the Input Space for Molecular Design}
\label{subsection:input-constraints-simple}

With the purpose of CAMD in mind, we next describe the input space constraints, $A, x \in \Omega$.  These constraints limit the search space to include structures which try to emulate physically-feasible molecular structures.  

\subsubsection{Basic MILP Formulation of Molecules}
\label{subsection:input-constraints-simple-basic}
QSPR methods used for property prediction in previous works are mostly based on group contribution methods \cite{zhang2015generic}.  As a result MILP formulations of molecules used in CAMD are also often based on group contribution methods.  
Modelling chemical properties with GNNs means that molecules are described in terms of an adjacency matrix $A \in \{0,1\}^{N \times N}$ and feature vectors $X \in \{0,1\}^{N \times F}$.  We therefore introduce an MILP formulation for molecules based on topological structure similar to the input of GNNs. 
This is similar in spirit to topological indexing methods for QSPR~\citep{austin2016computer}, which are in turn based on chemical graph theory. 

The structure of a solution is described by the adjacency matrix $A$, where $A_{ij} = 1$ indicates that node $i$ is connected to node $j$. The entries of a feature vector of a node $i$ are indicated by $x_{if}$, where $f$ is the position of a feature in that vector. The simplest machine learning model we consider comprises 14 features, which represent the knowledge summarised in Table~\ref{tab:features}.
\begin{table}[tb]
\caption{Background knowledge about molecule properties}
\begin{center}
\begin{tabular}{c c c} 
\toprule
 $X_{if}$ & type & descriptor \\
\midrule
 1 & atom & C \\
 2 & atom & O \\
 3 & atom & F \\
 4 & atom & Cl \\
 5 & neighbours & 0 \\
 6 & neighbours & 1 \\
 7 & neighbours & 2 \\
 8 & neighbours & 3 \\
 9 & neighbours & 4 \\
 10 & hydrogen & 0 \\ 
 11 & hydrogen & 1 \\ 
 12 & hydrogen & 2\\  
 13 & hydrogen & 3\\  
 14 & hydrogen & 4\\    
\bottomrule
\end{tabular}
\end{center}
\label{tab:features}
\end{table}%

With the adjacency matrix and the feature vectors for all the nodes, we introduce 
constraints to avoid trivially infeasible molecules. These constraints are summarised here and detailed in~\ref{appendix:molecule-constraints}:
\begin{enumerate}[label=(\alph*)]
    \itemsep0em 
    \item Molecules should be connected and of at least length 2. 
    \item Nodes are active if and only if they are connected to others. 
    \item To avoid redundancies, no gaps should exist between activated atoms (a molecule of length 3 should be $A_{11} = A_{22} = A_{33} = 1$ and not $A_{11} = A_{22} = A_{55} = 1$). 
    \item Each node should only have one atom type.
    \item The covalence of each atom must equal the number of active neighbours. 
\end{enumerate}

The formulation is not a tight formulation, and molecules can be found in the search space that might not be able to be synthesized or stable in a natural setting.  For instance, the formulation does not consider steric constraints on the bonds.  There are also molecules which are excluded from the search space.  An example of these are molecules with double or triple bonds.  Therefore we next describe how to extend this model. 

\subsubsection{Additional Properties}
The formulation in the previous subsection is seen as a basis which can be extended or modified for a particular CAMD setting. 
First of all, to limit the search space to molecules with no loops, consider the following constraint:
\begin{equation}
    \label{eq:numNodesMinusOne}
    \sum_i^{n-1} \sum_{j | j > i}^{n} A_{ij} = n - 1
\end{equation}

The constraint guarantees that the number of edges (LHS) is equal to the number of nodes minus one.  
\ifdefined\isarxiv
\else
When added to the set of constraints \eqref{eq:input-constraints}, no loops will be present in the molecules in the search space. Before adding the extra constraint, the search space only includes connected graphs due to constraint \eqref{eq:input-space-controversial-constraint}. 
\fi
This fact, combined with basic graph properties and the result of constraint \eqref{eq:numNodesMinusOne}, guarantees the graph to be acyclic.

Next, we consider constraints to have the search space include double bonded molecules. To achieve double bonds in our formulation, an extra feature is included in the feature vectors $X \in \mathbb{R}^{N \times F}$.  This feature $x_{i,15}$ indicates whether node $i$ is included in at least one double bond. This is a learnable parameter for the GNN. Outside the context of MILP formulations for GNNs, this feature would be included in a GNN which also includes edge features.  However, we leave such network architectures for subsequent work.

We introduce a binary variable $db_{il}$ that tracks whether a double bond is present between nodes $i$ and $l$. 
\begin{subequations}
\label{eq:double-bond-MILP}
\allowdisplaybreaks
    \begin{align}
        \label{eq:double-bond-MILP-0}
        &3 \cdot db_{il}                     \leq   x_{i,15} + x_{l,15} + A_{il}        & \forall i,l \in \{1, \dots, n\} \\
        \label{eq:double-bond-MILP-1}
        &2 \cdot x_{i,1} + 1 \cdot x_{i,2}       \geq   \sum_l^n db_{il}                    & \forall i \in \{1, \dots, n\} \\
        \label{eq:double-bond-MILP-2}
        &4 \cdot x_{i,1} + 2 \cdot x_{i,2} + 1 \cdot x_{i,3} + 1 \cdot x_{i,4} = \\
        & \sum_{s=0}^4 s \cdot x_{i,(5+s)} + \sum_{s=0}^4 s \cdot x_{i,(10+s)} + \sum_l^n db_{il} 
            & \forall i \in \{1, \dots, n\} \\
        \label{eq:double-bond-MILP-3}
        &x_{i,15}                        \leq \sum_l^n db_{il}                      & \forall i \in \{1, \dots, n\} \\
        \label{eq:double-bond-MILP-4}
        &db_{il}                         = db_{l,i}                                  & \forall i,l \in \{1, \dots, n\} \\
        \label{eq:double-bond-MILP-5}
        &db_{i,i}                         = 0                                        & \forall i \in \{1, \dots, n\}
    \end{align}
\end{subequations}

Constraint~\eqref{eq:double-bond-MILP-0} enforces that double bonds are only possible if nodes $i$ and $l$ are connected and $x_{i,15} = x_{l,15}= 1$. 
Constraint \eqref{eq:double-bond-MILP-1} limits the number of double bonds based on the covalence of the molecule. For instance, $x_{i,1}$ indicates a carbon molecule, meaning that there can be a maximum of $2$ double bonds. Constraint \eqref{eq:double-bond-MILP-2} limits the total number of bonds (including double bonds). 
Constraint \eqref{eq:double-bond-MILP-3} forces $x_{i,15}$ to be zero if there are no double bonds connected to node $i$.
The final two constraints are a symmetry constraint and a constraint indicating that a node cannot have a double bond with itself.

For triple bonds, the above formulation would be nearly identical. Instead, the variable $db_{il}$ would be replaced by $tb_{il}$ in all constraints to indicate a triple bond between node $i$ and node $l$. Constraint \eqref{eq:double-bond-MILP-1} would have $1 * x_{i,1}$ on the left-hand side because generally, only a carbon molecule can have a triple bond. Finally, in constraint \eqref{eq:double-bond-MILP-2}, $\sum_l^n tb_{il}$ would include a scalar multiple of $2$, because every triple bond removes two binding opportunities. 

Including both triple and double bonds in the search space can be achieved with the following set of constraints. 
\begin{subequations}
\label{eq:triple-double-bond-MILP}
\allowdisplaybreaks
    \begin{align}
        \label{eq:triple-double-bond-MILP-0-a}
        &3 \cdot db_{il}                     \leq   x_{i,15} + x_{l,15} + A_{il}        & \forall i,l \in \{1, \dots, n\} \\
        \label{eq:triple-double-bond-MILP-0-b}
        &3 \cdot tb_{il}                     \leq   x_{i,16} + x_{l,16} + A_{il}        & \forall i,l \in \{1, \dots, n\} \\
        \label{eq:triple-double-bond-MILP-1-a}
        &2 \cdot x_{i,1} + 1 \cdot x_{i,2}       \geq   \sum_l^n db_{il}                    & \forall i \in \{1, \dots, n\} \\
        \label{eq:triple-double-bond-MILP-1-b}
        &1 \cdot x_{i,1}                     \geq   \sum_l^n tb_{il}                    & \forall i \in \{1, \dots, n\} \\
        \label{eq:triple-double-bond-MILP-2}
        \begin{split}
            &4 \cdot x_{i,1} + 2 \cdot x_{i,2} + 1 \cdot x_{i,3} + 1 \cdot x_{i,4} = \\
            &\qquad \sum_{s=0}^4 s \cdot x_{i,(5+s)} + \sum_{s=0}^4 s \cdot x_{i,(10+s)} \\
            &\qquad + \sum_l^n db_{il} + \sum_l^n 2 \cdot tb_{il} 
        \end{split}
                                                                                    & \forall i \in \{1, \dots, n\} \\
        \label{eq:triple-double-bond-MILP-3-a}
        &x_{i,15}                        \leq \sum_l^n db_{il}                      & \forall i \in \{1, \dots, n\} \\
        \label{eq:triple-double-bond-MILP-3-b}
        &x_{i,16}                        \leq \sum_l^n tb_{il}                      & \forall i \in \{1, \dots, n\} \\
        \label{eq:triple-double-bond-MILP-4-a}
        &db_{il}                         = db_{li}                                  & \forall i,l \in \{1, \dots, n\} \\
        \label{eq:triple-double-bond-MILP-4-b}
        &tb_{il}                         = tb_{li}                                  & \forall i,l \in \{1, \dots, n\} \\
        \label{eq:triple-double-bond-MILP-5}
        &db_{ii}                         = tb_{ii} = 0                              & \forall i \in \{1, \dots, n\} \\
        \label{eq:triple-double-bond-MILP-6}
        &db_{li} + tb_{li}               \leq 1                                         & \forall i,l \in \{1, \dots, n\}
    \end{align}
\end{subequations}
where $x_{i,15}$ and $x_{i,16}$ indicate that an atom $i$ is part of a double or triple bond respectively. 

Once again we note that the introduction of these constraints does not span the entire space of possible molecules, nor does it include only naturally feasible molecules. For instance, introducing the triple bonds constraints would not find the molecule carbon monoxide.

\subsection{Bound Tightening Techniques}
\label{sec:method:fbbt}
Solving times of linear programming solvers are influenced by the tightness of the big-M constraints. It is therefore important to find tight constraints of the big-M values associated with a neuron.
We first take a look at the computationally efficient method of feasibility based bound tightening (FBBT). We first consider this for regular MLPs and then we continue adapting these methods for the GCN and GraphSAGE. 

\subsubsection{Big-M Coefficients for MLPs}
\label{subsubsection:fbbt-ann}
Feasibility based bound tightening techniques are bound tightening techniques which limit the feasible solution space by propagating the domain of the input space through the non-linear expression. This technique relies on interval arithmetic to compute the bounds on constraint activations over the variable domains \citep{gleixner2017three}. For the formulation of MLPs in Eq.~\eqref{eq:bigM}, recall we require bounds such that $l \leq w^Ty + b \leq u$. 
We now denote these as $l_j^k$ and $u_j^k$ for a node $j$ in layer $k$. 
Using interval arithmetic we can find bounds for the nodes for $k \geq 2$ in two ways, which result in the same bounds. For the first method, for layers $k \geq 2$ we find the upper bound $u_j^k$ and lower bound $l_j^k$ as follows:
\begin{subequations}
\begin{align}
\label{eq:find-upper-bound-ANN}
u_{j}^{k} &=\sum_{i=1}^{n_{k-1}} \max \left\{w_{j i}^{k} \max \left\{0, u_{i}^{k-1}\right\}, w_{j i}^{k} \max \left\{0, l_{i}^{k-1}\right\}\right\}+b_{j}^{k}, \\
l_{j}^{k} &=\sum_{i=1}^{n_{k-1}} \min \left\{w_{j i}^{k} \max \left\{0, u_{i}^{k-1}\right\}, w_{j i}^{k} \max \left\{0, l_{i}^{k-1}\right\}\right\}+b_{j}^{k}.
\end{align}    
\end{subequations}

Note that the inner $\max$ operators function as the ReLU activation function. The outer $\max$ function is necessary since the weight matrix entries can also be negative. For $k = 1$ we remove the inner $\max$ functions as the input is not necessarily positive since there is no ReLU operator. The same bounds can be found by solving the LP problems:
\begin{equation}
\label{eq:feasibility-bounds-ann-a}
\begin{aligned}
&u_{j}^{k}=\max \left\{t_{j}^{k}: t_{j}^{k} \in C_{j}^{k}\right\} \\
&l_{j}^{k}=\min \left\{t_{j}^{k}: t_{j}^{k} \in C_{j}^{k}\right\}
\end{aligned}    
\end{equation}
for the constraint set
\begin{equation}
\label{eq:feasibility-bounds-ann-b}
\begin{aligned}
C_{j}^{k}=&\left\{t_{j}^{k}: t_{j}^{k}=w_{j}^{k} x^{k-1}+b^{k}, x^{k-1} \in\left[\max \left\{0, L^{k-1}\right\}, \max \left\{0, U^{k-1}\right\}\right] 
\subset \mathbb{R}^{n_{k-1}}\right\}
\end{aligned}    
\end{equation}

To speed up the solving time, some activation variables $z$ can be determined based on the value of the lower and upper bound.  When the lower bound $l_j^k$ of a particular node is above $0$, $z_j^k$ can be set to $1$. In this case, it is known that $x_j^k$ will always be positive and thus $z_j^k$ must be $1$.  The same goes for a positive lower bound $l_j^k > 0$: in this case $z_j^k = 0$. 

\subsubsection{Big-M Coefficients for GCNs}
The following subsection explains how to find the upper and lower bounds associated with the ReLU constraints for a GCN model. 
Specifically, we require the upper and lower bounds, $U_{ij}^k$ and $L_{ij}^k$ in Eq.s~\eqref{eq:conditional-hidden-layer}--\eqref{eq:GCN-MILP-c}.

It is assumed that the lower and upper bounds of all the input feature vectors are the same. This is because the input feature vectors of all nodes describe the same features of those nodes, and the feature vectors must be equal in length. This makes the bound propagation symmetric over all nodes, which in turn allows us to only calculate the bounds of all nodes once per layer $k$. 
Before the optimisation, for node $i$, the number of neighbouring nodes and the number of their respective neighbours are unknown. In the case of maximisation, we have to find a scalar which upper bounds $H_{ij}^k$ for all possible neighbourhood structures of node $i$. Specifically, we compute $d_i^+$ and $d_l^+$ such that the following is maximised (remember that $d_i^+ = d_i + 1$, where $d_i$ is the degree of node $i$ if self loops are not possible):

\begin{equation}
\begin{aligned}
    & \underset{d_i^+, d_l^+}{\max} &d_i^+ \frac{1}{\sqrt{d^{+}_i d^{+}_l}} {\bm{W}^{k}_j}^T \bm{H}_l^{(k-1)}
\end{aligned}
\end{equation}

The upper bound of ${\bm{W}^{k}_j}^T \bm{H}_l^{(k-1)}$ is determined as in  Eq.~\eqref{eq:find-upper-bound-ANN} with zero bias. If this upper bound is positive, we want to add as much as possible to account for all possible neighbourhood structures, i.e., $\frac{d_i^+}{d_l^+}$ should be maximized. 
The degree of node $i$, $d_i^+$, can maximally be $d_{max} + 1$, and the minimum degree of the neighbours of $i$ needs to be $d_{l}^+ = 1 + 1$. 
One of other hand, if the upper bound from Eq.~\eqref{eq:find-upper-bound-ANN} is negative, $\frac{d_i^+}{d_l^+}$ should be minimized following the same logic. This results in 
\begin{subequations}
\begin{align}
& \begin{split}
 U_{ij}^{k} = \max & \bigg\{ \sqrt{\frac{{d_{max} + 1}}{{2}}} \sum_{s=1}^{n_{k-1}} \max \left\{w_{j s}^{k} \max \left\{0, U_{sj}^{k-1}\right\}, w_{j s}^{k} \max \left\{0, L_{sj}^{k-1}\right\}\right\}, 
 \\ & \sqrt{\frac{{2}}{{d_{max} + 1}}} \sum_{s=1}^{n_{k-1}} \max \left\{w_{j s}^{k} \max \left\{0, U_{sj}^{k-1}\right\}, w_{j s}^{k} \max \left\{0, L_{sj}^{k-1}\right\}\right\} \bigg\}
\end{split} \\
& \begin{split}
 L_{ij}^{k} = \min & \bigg\{ \sqrt{\frac{{d_{max} + 1}}{{2}}} \sum_{s=1}^{n_{k-1}} \min \left\{w_{j s}^{k} \max \left\{0, U_{sj}^{k-1}\right\}, w_{j s}^{k} \max \left\{0, L_{sj}^{k-1}\right\}\right\}, 
 \\ & \sqrt{\frac{{2}}{{d_{max} + 1}}} \sum_{s=1}^{n_{k-1}} \min \left\{w_{j s}^{k} \max \left\{0, U_{sj}^{k-1}\right\}, w_{j s}^{k} \max \left\{0, L_{sj}^{k-1}\right\}\right\} \bigg\}
\end{split}
\end{align}
\end{subequations}

In practice, we only need the first term of each bound, as in the case that $\max \left\{w_{j i}^{k} \max \left\{0, u_{i}^{k-1}\right\}, w_{j i}^{k} \max \left\{0, l_{i}^{k-1}\right\}\right\}$ is negative, the node is turned off and on for the upper bound and lower bound respectively (as described in Section~\ref{subsubsection:fbbt-ann}).

For a similar formulation as Eqs.~\eqref{eq:feasibility-bounds-ann-a} and \eqref{eq:feasibility-bounds-ann-b}, the following set of equations can be considered:
\begin{subequations}
\begin{align}
    U_{ij}^{k}&=\max \left\{t_{ij}^{k}: t_{ij}^{k} \in C_{ij}^{k}\right\} \\
    L_{ij}^{k}&=\min \left\{t_{ij}^{k}: t_{ij}^{k} \in C_{ij}^{k}\right\} \\
    C_{ij}^{k}&=\left\{t_{ij}^{k}: t_{ij}^{k}= \sqrt{\frac{{d_{max} + 1}}{{2}}} w_{j}^{k} x^{k-1}, \right.\\
        &\qquad\left. x^{k-1} \in\left[\max \left\{0, L_{i}^{k-1}\right\}, \max \left\{0, U_{i}^{k-1}\right\}\right]\subset \mathbb{R}^{n_{k-1}}\right\}
\end{align}
\end{subequations}
 
\subsubsection{GraphSAGE}
\label{subsubsubsection:FBBT-GraphSAGE}
For GraphSAGE, FBBT is very similar to the FBBT proposed for the GCN. 
We now require the upper and lower bounds, $U_{ij}^k$ and $L_{ij}^k$ in Eq.s~\eqref{eq:graph-conv-relu-1}--\eqref{eq:graph-conv-relu-3}.

There are two parts that contribute to the total bound. The first, which is associated with the root node $i$, is calculated similarly as the MLP FBBT. The second part is calculated in the same way as the GCN, but with the root node $i$ omitted. This results in the following bound propagation equations:
\begin{subequations}
\begin{align}
 \begin{split}
     U_{ij}^{k} &= \sum_{s=1}^{n_{k-1}} \max \left\{\hat{w}_{j s}^{k} \max \left\{0, U_{sj}^{k-1}\right\}, \hat{w}_{j s}^{k} \max \left\{0, L_{sj}^{k-1}\right\}\right\} \\&+ \sqrt{\frac{{d_{max}}}{{2}}} \sum_{s=1}^{n_{k-1}} \max \left\{\bar{w}_{j s}^{k} \max \left\{0, U_{sj}^{k-1}\right\}, \bar{w}_{j s}^{k} \max \left\{0, L_{sj}^{k-1}\right\}\right\},
 \end{split} \\
 \begin{split}
     L_{ij}^{k} &= \sum_{s=1}^{n_{k-1}} \min \left\{\hat{w}_{j s}^{k} \max \left\{0, U_{sj}^{k-1}\right\}, \hat{w}_{j s}^{k} \max \left\{0, L_{sj}^{k-1}\right\}\right\} \\&+ \sqrt{\frac{{d_{max}}}{{2}}} \sum_{s=1}^{n_{k-1}} \min \left\{\bar{w}_{j s}^{k} \max \left\{0, U_{sj}^{k-1}\right\}, \bar{w}_{j s}^{k} \max \left\{0, L_{sj}^{k-1}\right\}\right\}.
 \end{split}
\end{align}
\end{subequations}
Notice that $\sqrt{\frac{{d_{max}}}{{2}}}$ replaces $\sqrt{\frac{{d_{max}}}{{2}}}$ because the root node is omitted. 

For a similar formulation as Eqs.~\eqref{eq:feasibility-bounds-ann-a} and \eqref{eq:feasibility-bounds-ann-b}, the following set of equations can be considered:
\begin{subequations}
\begin{align}
    U_{ij}^{k}&=\max \left\{t_{ij}^{k}: t_{ij}^{k} \in C_{ij}^{k}\right\} \\
    L_{ij}^{k}&=\min \left\{t_{ij}^{k}: t_{ij}^{k} \in C_{ij}^{k}\right\} \\
    C_{ij}^{k}&=\left\{t_{ij}^{k}: t_{ij}^{k}= \hat{w}_{j}^{k} y^{k-1} + \sqrt{\frac{{d_{max}}}{{2}}} \bar{w}_{j}^{k} x^{k-1}, \right.\\
        &\qquad\left. x^{k-1},y^{k-1} \in\left[\max \left\{0, L_{i}^{k-1}\right\}, \max \left\{0, U_{i}^{k-1}\right\}\right]\subset \mathbb{R}^{n_{k-1}}\right\}
\end{align}
\end{subequations}
where $\hat{w}$ is the weight matrix associated with the root node, and $\bar{w}$ is the weight matrix multiplied with the aggregated nodes. 

The upper bounds calculated in this subsection also serve as the values for the big-M constraints in \eqref{eq:graph-conv-conditional-support-1}, \eqref{eq:graph-conv-conditional-support-2} and \eqref{eq:graph-conv-conditional-support-3}.  For instance, consider constraint \eqref{eq:graph-conv-conditional-support-1} and let there be no connection between node $i$ and node $l$. In that case $b_{il}^k$ needs to be able to be equal to 0. Enough $\bm{M}$ needs to be subtracted such that the left hand side of the equation is lower than 0. To achieve this $\bm{M}$ needs to be larger than $H_l^k$, which can be achieved if $\bm{M}$ is equal to the upper bound of $H_l^k$.

\subsection{Benchmark Genetic Algorithm}
As a benchmark, we implement a straightforward genetic algorithm (GA) to optimise over trained GNNs.  This GA does not depend on a latent space architecture, as proposed by \citet{rittig2022graph}.  It uses a string representation of the symmetric adjacency matrix of the molecules and of the feature vectors of the molecules such that single-point crossovers and string mutations can be applied.  A description is given in \ref{appendix:ga}.

\section{Numerical Results}
\label{section:results}\label{sec:results}

Recall that the goal of this work is to formulate graph neural networks such that they can be used as surrogate models in optimisation problems.  To validate the MI(N)LP formulations, we turn to the case study of chemical property prediction, specifically the prediction of boiling points.

\subsection{Experimental Setup}
\label{subsection:results-experimental-setup}
The workflow of the experiments is as follows. First a data set is chosen, which will be used to train a GNN. Thereafter, this trained GNN is represented using the methods described in Section~\ref{section:methods}. The resulting formulations are optimised using a deterministic solver, namely Gurobi version 9.5.1 \citep{gurobi-951}.

The experiments were ran on two different machines, a laptop, and a virtual machine. The laptop was used for quick, low-resource-intensive experiments in which we were only interested in the MI(N)LP solutions. 
The virtual machine was used for experiments where solving times were compared. These experiments took longer and required constant CPU availability. The laptop was equipped 
with a 1.4 GHz Quad-Core Intel Core i5 processor, 8GB 
memory. 
The virtual machine was equipped with eight 2.5 Ghz Intel Xeon Gold 6248 CPUs, and 16GB memory. 
The machine learning models were trained using PyTorch 1.11.0 \citep{NEURIPS2019_9015}, and implemented using PyTorchGeometric 2.0.4 \citep{fey2019fast}.  

\subsubsection{Case Study}
\label{Section:Case-study-experimental setup}

In order to examine the output results of the proposed methods, we consider a representative test case where a chemist models a chemical property and tries to optimise it. This is a common use case of CAMD \cite{fruhbeis1987computer}. The chemist can then use the found solutions and test the predicted properties instead of having to search over the entire search space of feasible molecules. 
Models and training hyperparameters were selected heuristically following some preliminary experiments discussed in \ref{appendix:initial-experiments}. 
Note that accurately training GNNs is not the focus of this work. 

The input space constraints were extended by including the possibility to find double and triple bonds, implemented using the input constraints described by Eqs.~\eqref{eq:triple-double-bond-MILP}. Initial testing of the experiments found problems with steric constraints on the model, and to circumvent this, molecules with loops were excluded from the search space. 

The case study included the optimisation of the simplest GraphSAGE configuration with at least one hidden layer, which was the $\verb|1|$ $\times$
$\verb|16|$ configuration.  A total of three different formulations were tested: (i) only single and double bonds allowed, (ii) only single and triple bonds allowed, and (iii) single, double, and triple bonds allowed. 
As double and triple bonds are handled using extra features, three different neural neural network configurations were trained.
The GNNs are trained using the best hyperparameters found in our preliminary experiments, and the model with the lowest validation error was selected as input for the MILP formulation. The formulation was solved to optimality on the laptop. Multiple solutions (max 8) were recorded for molecule lengths $4 \text{ and } 5$. All solutions were analysed, checking whether the found solutions were molecules which naturally exist in nature. 

\subsection{Results}
The following section describes the results that were found in the experiments. First, the results are discussed for the GCN and GraphSAGE models individually. The results per model differ in node depth, layer width and molecule length. Thereafter, a comparison is described between models with the same parameters. Lastly, the CAMD case study is described, showing which molecules were found using the proposed methods.

\subsubsection{Initial Experiments}
\label{subsubsection:results-initial-tests}
\begin{figure}[tb!]
    \centering
    \includegraphics[clip,trim={2.4cm 0 0.2cm 2.6cm},width=0.8\textwidth]{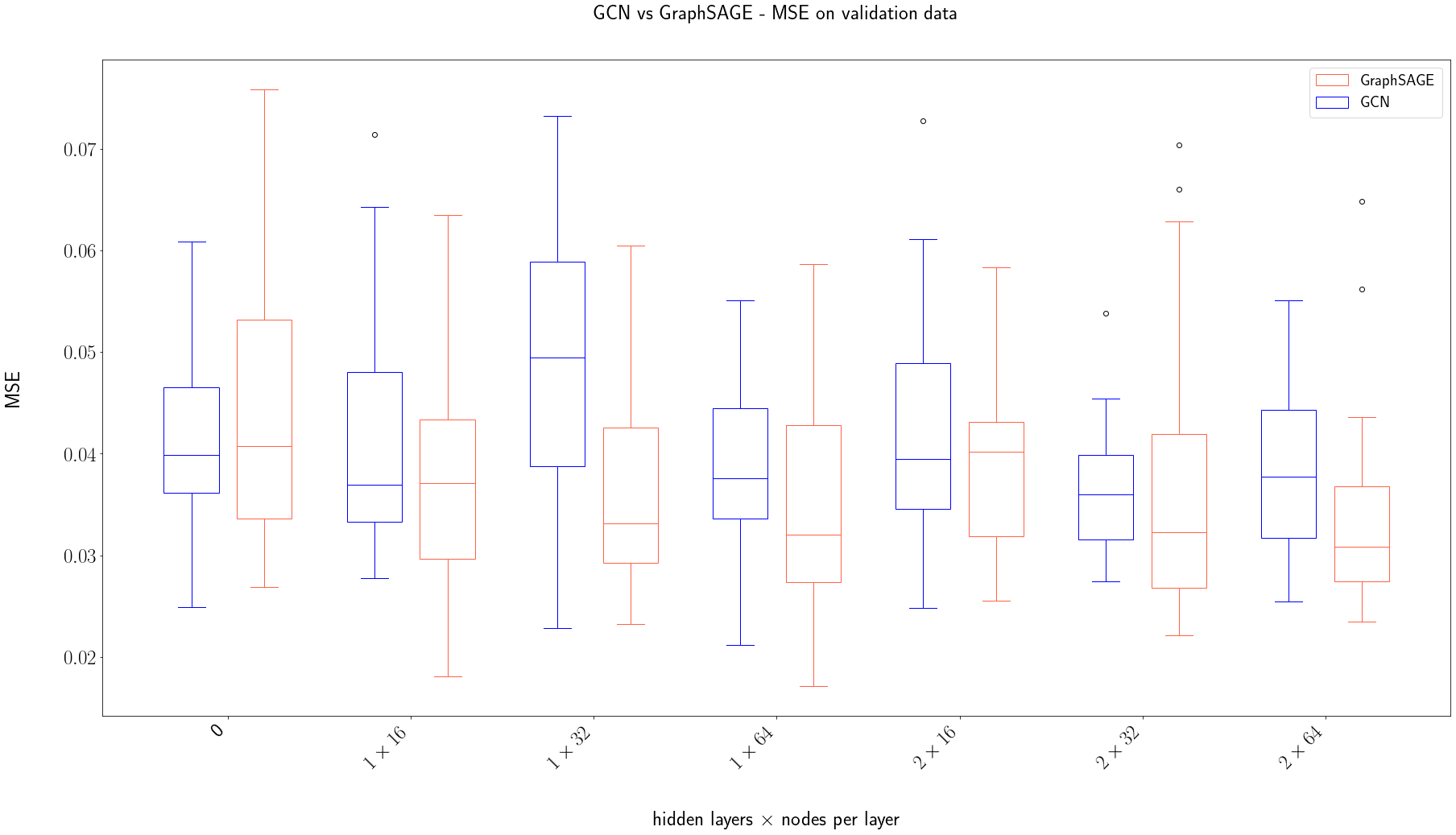}
    \caption{Box-plots of the MSE, for the GCN and GraphSAGE models for different layer depths and node width, after independently running the models 20 times for each configuration. The box-plots indicate the median, the lower and upper quartile and the lowest and largest MSE, when outliers, which are indicated by the dots, are excluded.}
    \label{fig:boxplots}
\end{figure}

For both GCN and GraphSAGE models, we trained seven configurations (differing in the amount of hidden layers and nodes per hidden layer) using the MSE as a loss function.  A comparison of the box plots of the MSE of the models can be found in Fig.~\ref{fig:boxplots}. 

For the GCN model, the minimum MSE values range from $0.0213$ to $0.0278$ and the median MSE values range from $0.0360$ to $0.0495$.  For the GraphSAGE model the minimum MSE values range from $0.0171$ to $0.0277$ and the median from $0.0308$ to $0.0438$.  In case of the GraphSAGE we can see that median MSE values decrease as the amount of nodes per layer increase. For GCN no such pattern can be detected. 
Comparing the model accuracy of both GNN models show that the GraphSAGE model has a better median MSE validation value for 4 out of the 7 configurations. The minimum MSE values also show the GraphSAGE model to be better for 4 out of the 7 configurations. The best overall minimum MSE is found in the \verb|1 x 64| configuration of the GraphSAGE model, with an MSE of $0.0171$. However, the \verb|1 x 16| is a close second with an MSE of $0.0181$. 

\begin{table}[tp!]
    \caption{Solving times (up to 36000s; smaller is better) and remaining optimality gaps (smaller is better) for the GCN MINLP formulation for different molecule lengths.}
    \centering
    \footnotesize
    \begin{tabular}{c c c@{\hskip 4mm} ccc@{\hskip 4mm} ccc}
    \toprule
     & num. layers & 0 & \multicolumn{3}{c}{1} & \multicolumn{3}{c}{2} \\ \midrule
     & nodes/layer & 0 & 16 & 32 & 64 & 16 & 32 & 64 \\ 
     \midrule
     \multirow{2}{2cm}{molecule length 4} & Time (s) & 45 & 536 & 6426 & 13404 & 9428 & - & -\\
     & MIP Gap & 0.00 & 0.00 & 0.00 & 0.00 & 0.00 & 3.59 & 30.19\\ 
     \midrule
      \multirow{2}{2cm}{molecule length 6} & Time (s) & 2309 & - & - & - & - & - & -\\
     & MIP Gap & 0.00 & 2.92 & 8.61 & 8.20 & 31.61 & 47.10 & 93.86 \\ 
     \midrule
      \multirow{2}{2cm}{molecule length 8} & Time (s) & - & - & - & - & - & - & -\\
     & MIP Gap & 0.54 & 7.79 & 14.21 & 12.11 & 39.32 & 63.47 & 104.51 \\ 
     \bottomrule
    \end{tabular}
    \label{tab:results-GCN}
\end{table}

We then optimized the MILP and MINLP formulations of the trained GNNs with a maximum time limit of 10 hours.  Table~\ref{tab:results-GCN} compares the results for the optimisation of the MINLP formulation for the different configurations.  Only six out of the 21 experiments were solved within 10 hours. 
All configurations were solved for molecule length $n = 4$ apart from \verb|2 x 32| and \verb|2 x 64|, with solution times increasing with the number of nodes. For $n = 6$ only the formulation with $0$ hidden layers was solved to optimality. 
The optimality gaps are non-zero for the configurations where an optimum is not found. Recall that the optimality gap indicates how far the upper bound is removed from the lower bound, expressed in multiples of the lower bound. 
As the node depth increases, the remaining optimality gap after 10 hours increases, apart from increasing the node depth from $32$ to $64$ with $1$ hidden layer for $n = 6$ and $n = 8$. As the layer depth increases, the optimality gap also increases for all non-solved configurations. Finally, we note that all optimality gaps increase as the molecule length increases.

\begin{table}[tp!]
    \caption{Solving times (up to 36000s; smaller is better) and remaining optimality gaps (smaller is better) for the GraphSAGE MINLP formulation for different molecule lengths.}
    \centering
    \footnotesize
    \begin{tabular}{c c c@{\hskip 4mm} ccc@{\hskip 4mm} ccc}
    \toprule
     & num. layers & 0 & \multicolumn{3}{c}{1} & \multicolumn{3}{c}{2} \\ \midrule
     & nodes/layer & 0 & 16 & 32 & 64 & 16 & 32 & 64 \\
     \midrule
     \multirow{2}{2cm}{molecule length 4} & Time (s) & 3 & 384 & 1045 & 30248 & - & - & -\\
     & MIP Gap & 0.00 & 0.00 & 0.00 & 0.00 & 0.91 & 4.93 & 13.99 \\ 
     \midrule
      \multirow{2}{2cm}{molecule length 6} & Time (s) & 1108 & - & - & - & - & - & -\\
     & MIP Gap & 0.00 & 0.57 & 2.29 & 8.49 & 5.57 & 11.39 & 17.26 \\ 
     \midrule
      \multirow{2}{2cm}{molecule length 8} & Time (s) & - & - & - & - & - & - & -\\
     & MIP Gap & 0.33 & 5.19 & 3.63 & 11.18 & 7.11 & 13.60 & 26.54 \\ 
     \bottomrule
    \end{tabular}
    \label{tab:results-GraphSAGE}
\end{table}

Table~\ref{tab:results-GraphSAGE} shows the results of optimising the MILP formulations of the trained GraphSAGE neural networks. Only 5 configurations were solved to optimality. All the others terminated after a time limit of 10 hours. Of the solved cases, four optima were found when the search space was limited to atoms of length 4 ($n = 4$), the other one was found when $n = 6$. 
Once again, the solution times increase with both the size of the GNN and the length of the molecules.
We see that for $n = 4, 6$, as the node depth increases, the remaining optimality gaps become larger. For $n = 8$, this is not the case. When increasing the node depth for one hidden layer from $16$ to $32$, the optimality gap decreases. When increasing the number of layers, the optimality gap always gets larger when the node depth stays the same, for all molecule lengths. Finally as the molecule length increases, the optimality gap also increases.

We can also compare between the GCN and GraphSAGE models. As mentioned previously, six of the GCN formulations were solved within 10 hours, while only five of the GraphSAGE formulations were solved. 
For the solved problems, the GraphSAGE model is generally faster, apart from the configurations \verb|1 x 64| and \verb|2 x 16| for $n = 4$, where in the latter case GCN solves to optimality and GraphSAGE does not. 
We also obseve that GraphSAGE has a better optimality gap than the GCN formulation for most configurations and molecule lengths. There are a few exceptions however. For $n = 4$, \verb|2 x 16| and \verb|2 x 32|, GCN has a smaller optimality gap than the GraphSAGE formulation. The same goes for $n = 6$ with 1 hidden layer and 64 nodes in that layer. 

\begin{figure}[tbp]
    \centering
    \subfigure[GA vs GCN, n = 4]{\includegraphics[width=0.47\textwidth]{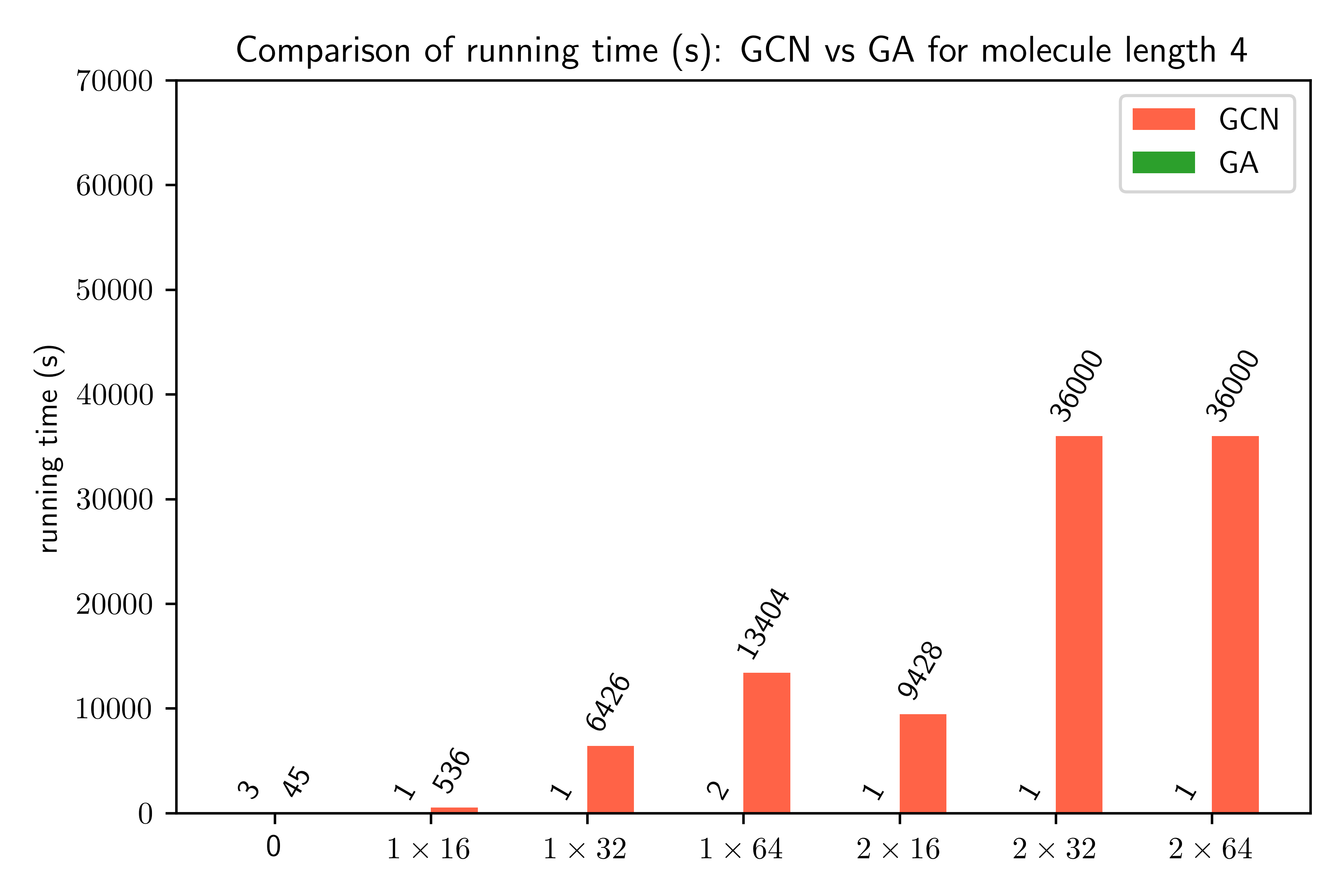}\label{fig:Ex_Im9}}\hspace{0.03\textwidth}
    \subfigure[GA vs GraphSAGE, n = 4]{\includegraphics[width=0.47\textwidth]{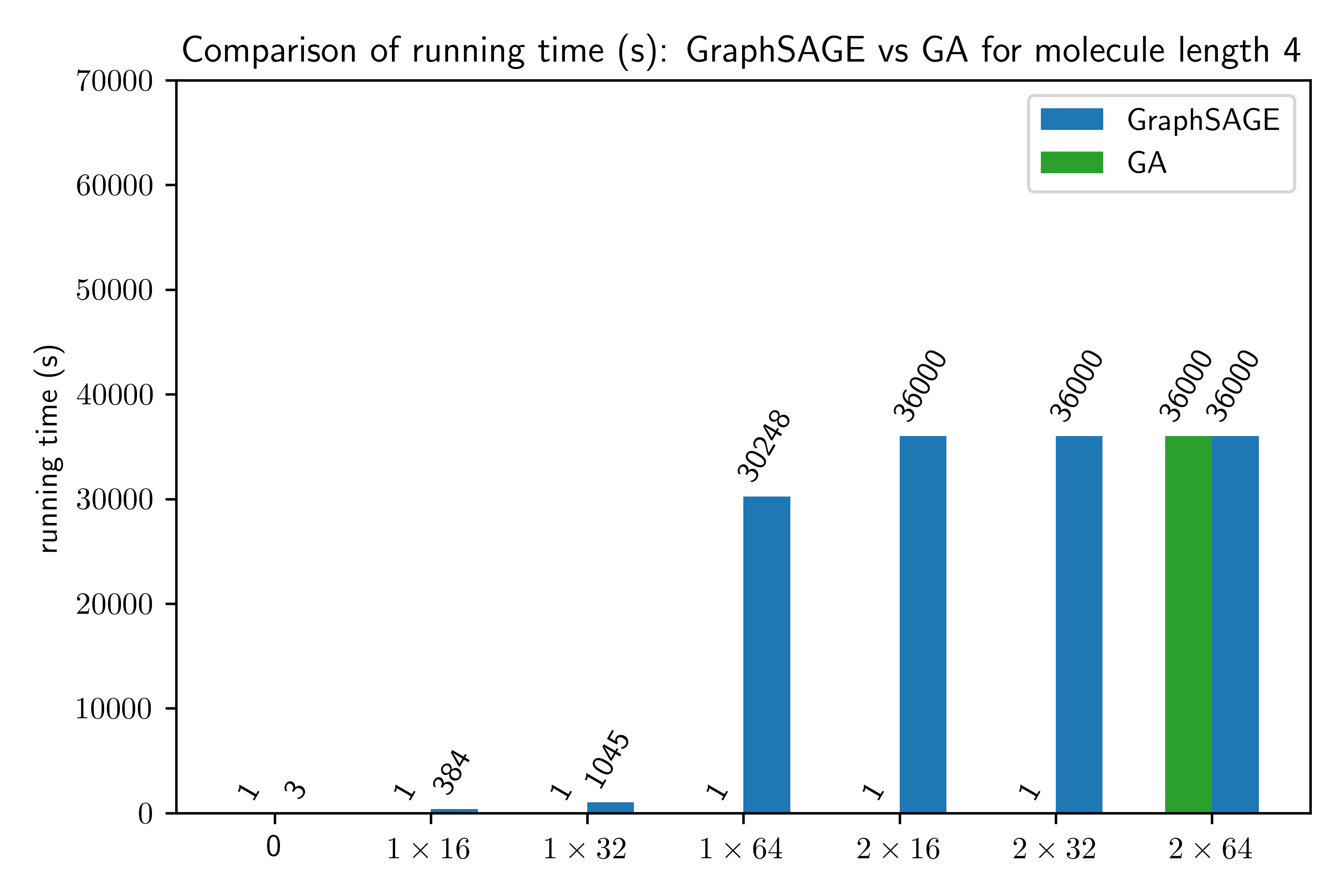}\label{fig:Ex_Im12}}
    
    \subfigure[GA vs GCN, n = 6]{\includegraphics[width=0.47\textwidth]{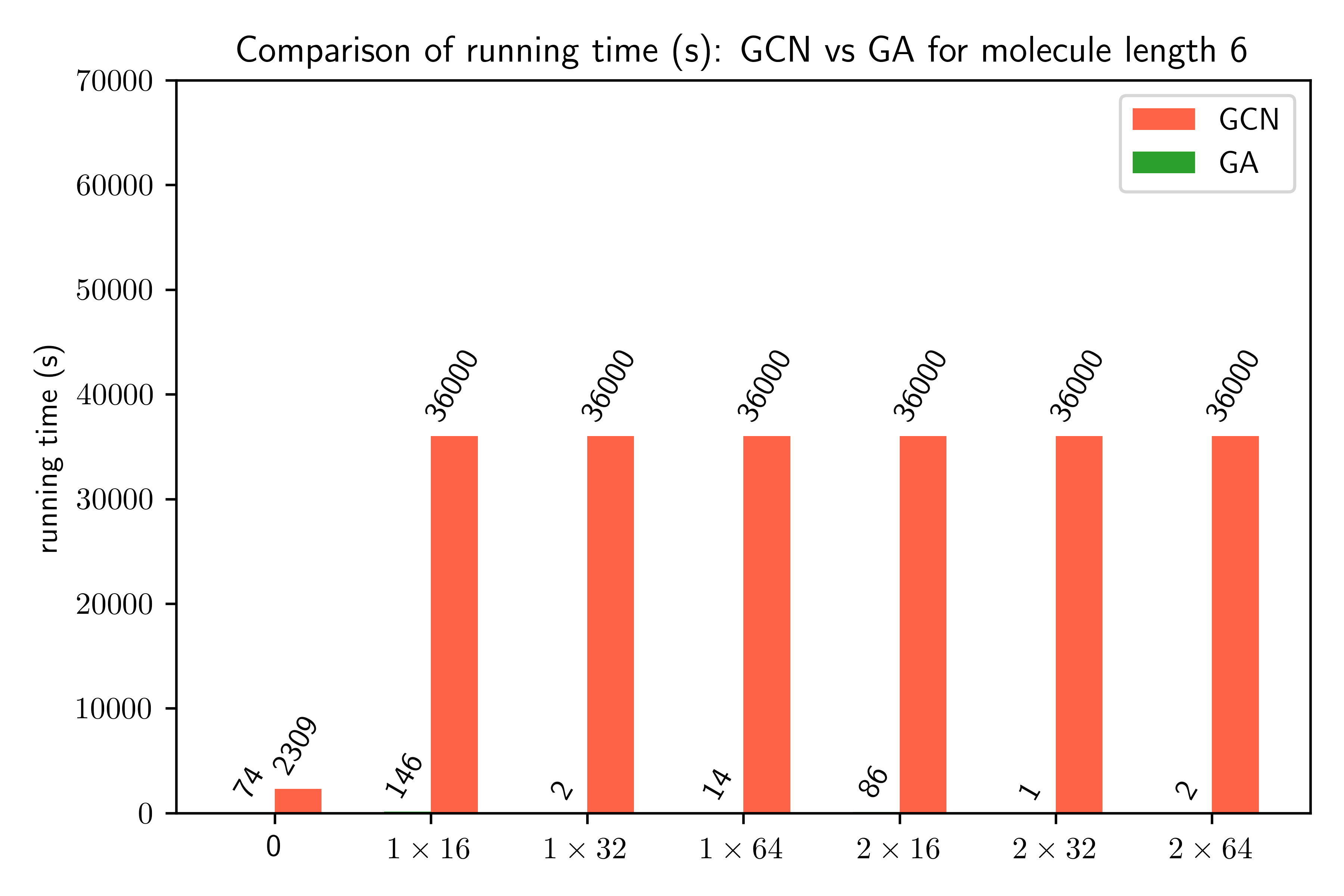}\label{fig:Ex_Im10}}\hspace{0.03\textwidth}
    \subfigure[GA vs GraphSAGE, n = 6]{\includegraphics[width=0.47\textwidth]{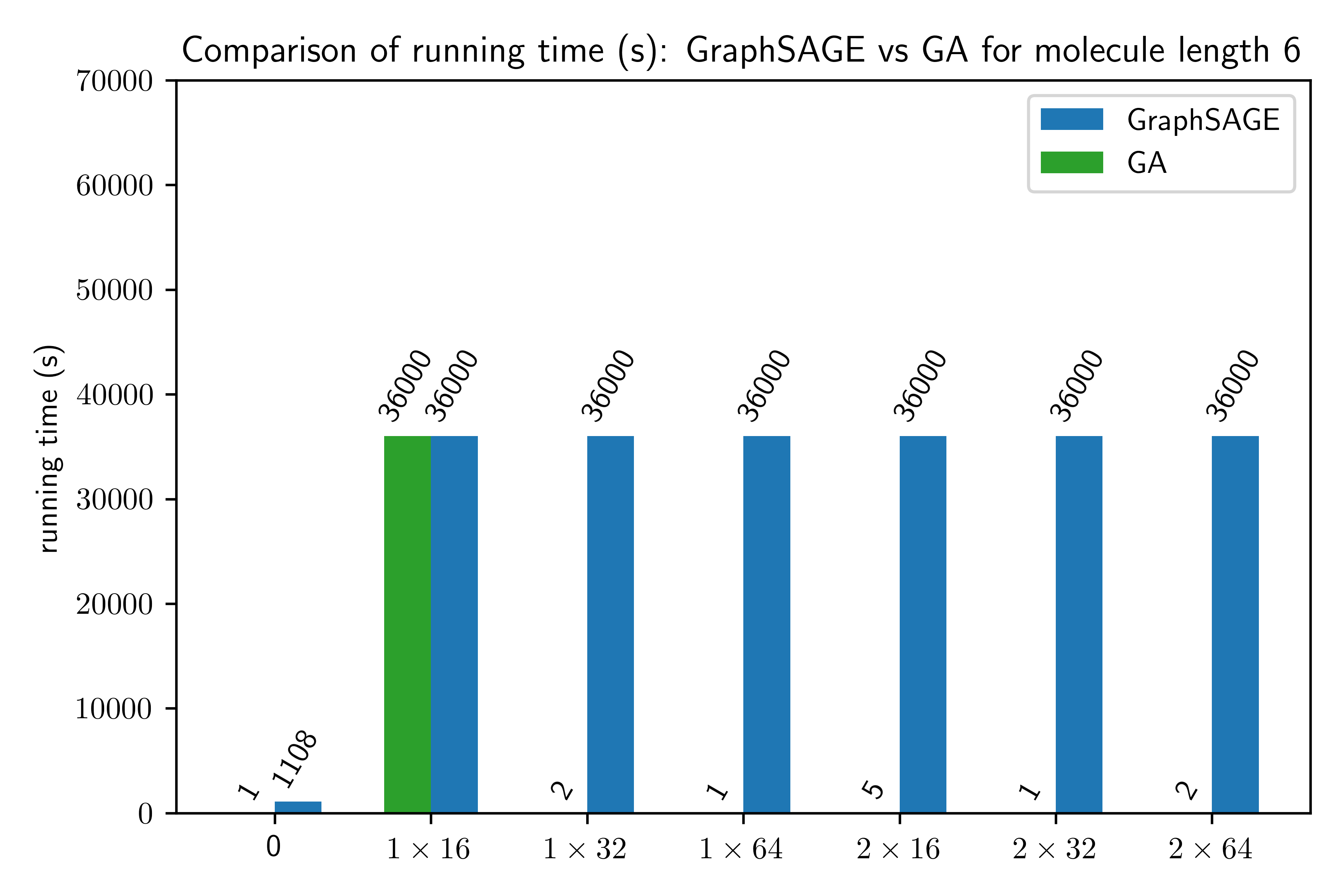}\label{fig:Ex_Im13}}
    
    \subfigure[GA vs GCN, n = 8]{\includegraphics[width=0.47\textwidth]{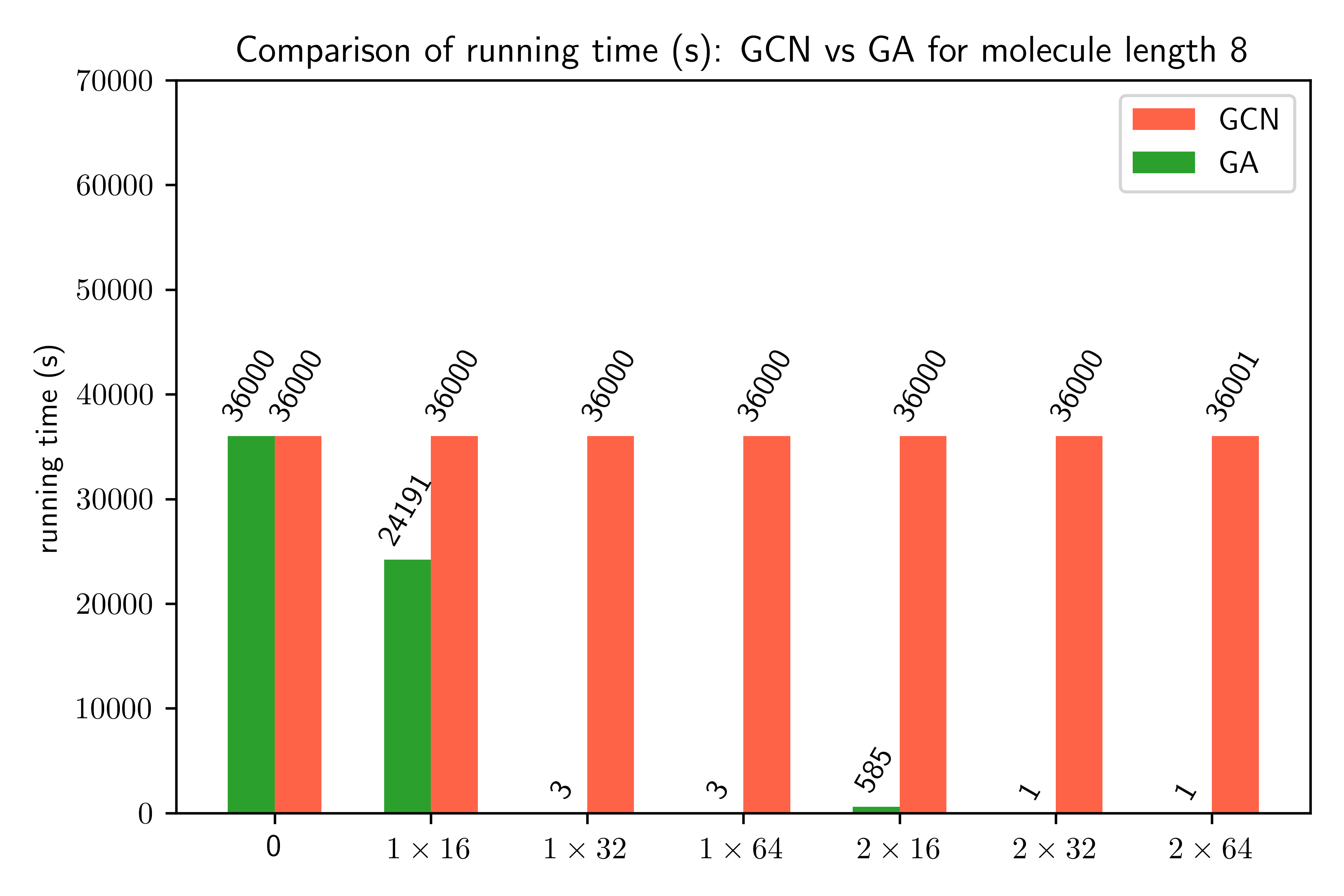}\label{fig:Ex_Im11}}\hspace{0.03\textwidth}
    \subfigure[GA vs GraphSAGE, n = 8]{\includegraphics[width=0.47\textwidth]{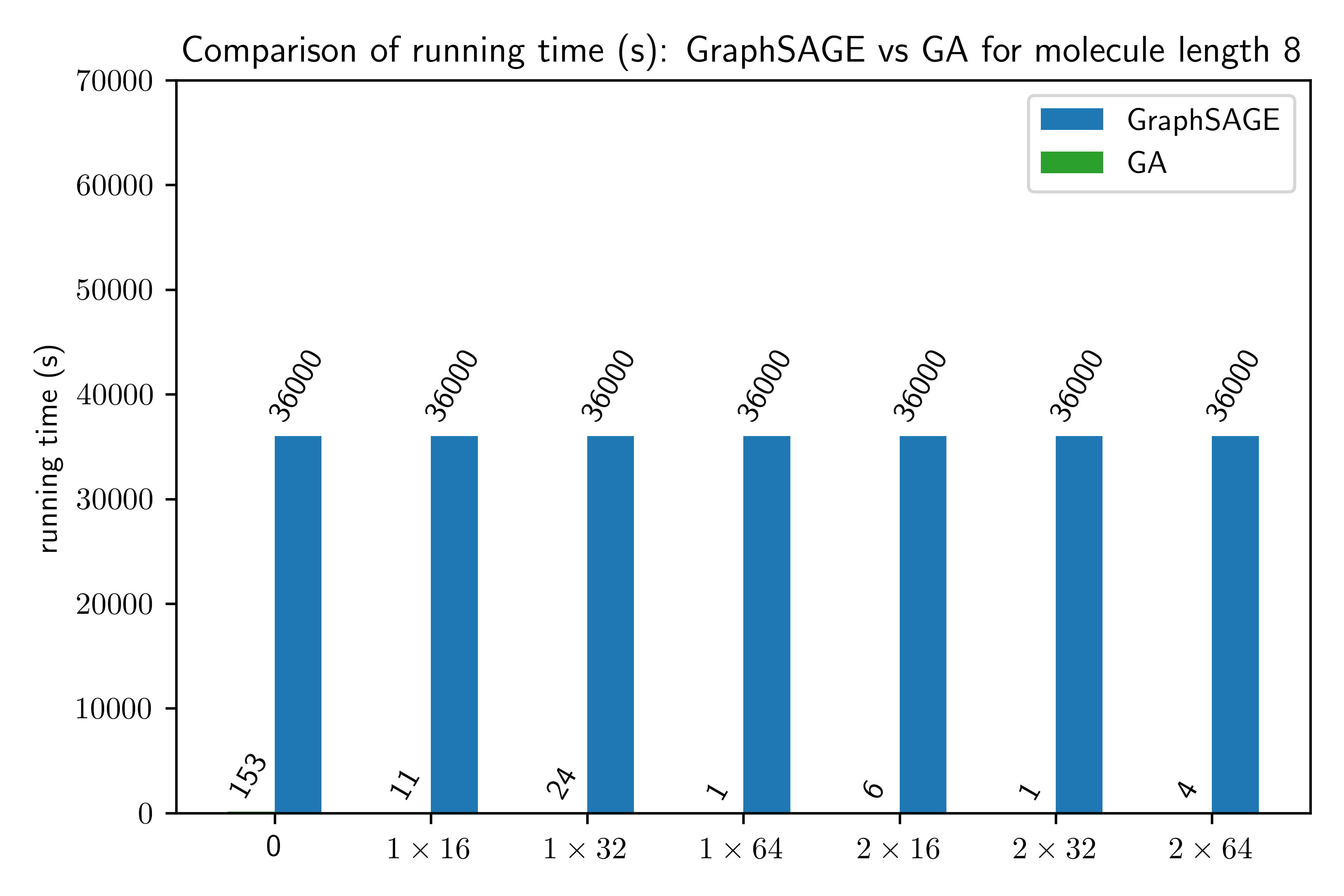}\label{fig:Ex_Im14}}    
    
    \caption{Comparison of the GA and GNN.  The solving time in seconds for the GA indicates after how many seconds the GA found an objective value of equal quality or better, than the MILP formulation of the GNN.}
    \label{fig:results-seconds-GA-vs-GNN}
\end{figure}

Finally we compare the solving times of the GNNs with a baseline.  Fig~\ref{fig:results-seconds-GA-vs-GNN} shows this comparison.  For 35 out of 42 instances ($83 \%$), the GA found an equally good solution as the deterministic optimiser in less than 2 minutes.  For 3 instances, the GA did not find an equally good objective value and terminated because it reached the time limit of $10$ hours. 

\subsubsection{Case Study}

\begin{figure}[p]
    \centering
    \hspace*{-0.5cm}\includegraphics[clip,trim={4.25cm 0 1.9cm 0},width=1.075\textwidth]{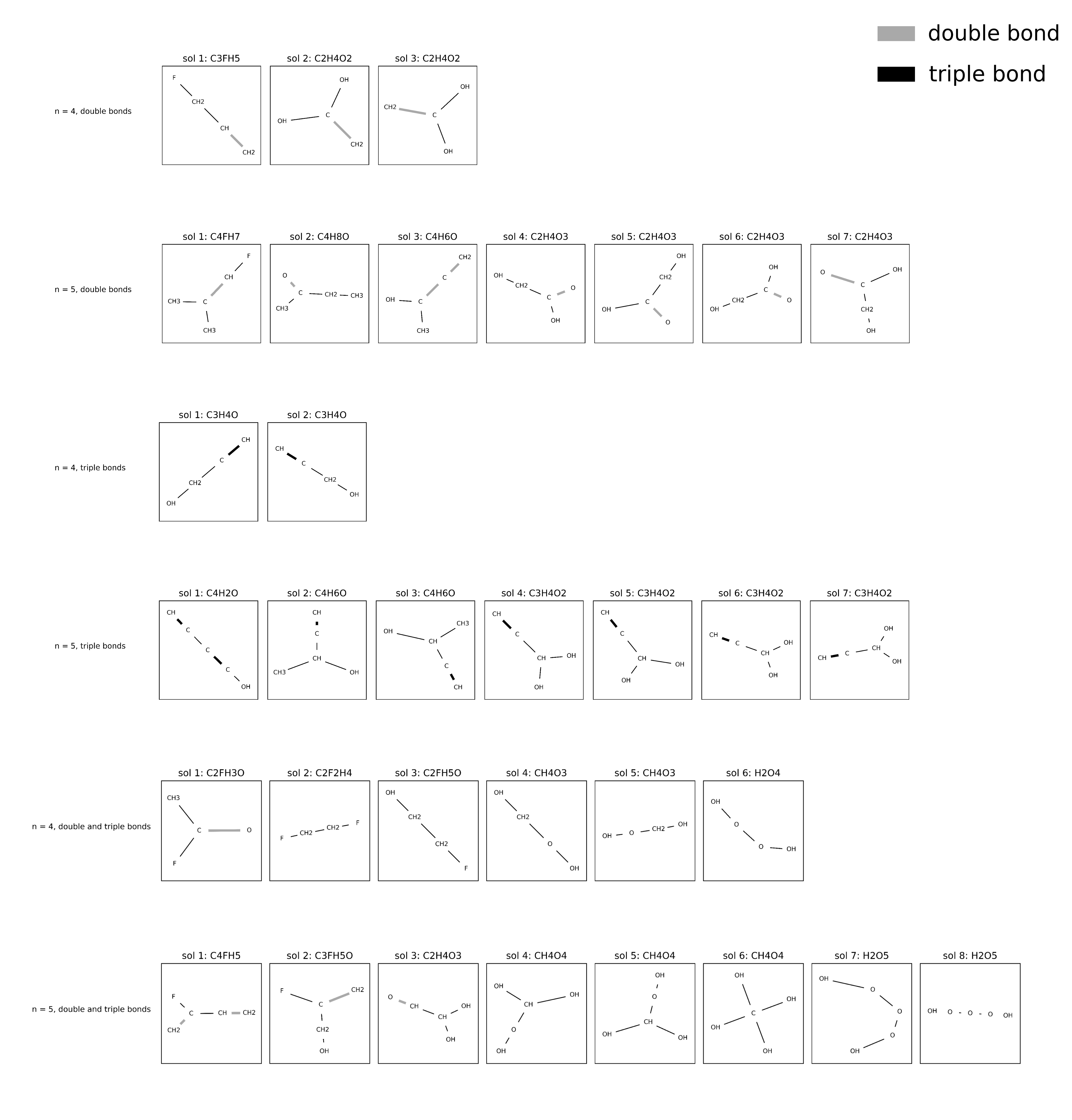}
    \caption{The molecules associated with the best found lower bounds during the branch and bound optimisation process of the GraphSAGE formulation with 1 hidden layer and 16 nodes.}
    \label{fig:result-molecules}
\end{figure}

For the case study we selected the GraphSAGE GNN with 1 hidden layer and 16 nodes.  The reason for choosing this model is explained in \ref{Appendix:extra-results-case-study-model-selection}.  Recall from Section~\ref{Section:Case-study-experimental setup} there were three different search spaces for the molecules.  These were the search space of single and double bonded molecules, single and triple bonded molecules, and single, double and triple bonded molecules.  Graphical representations of the found solutions for $n = 4$ and $n = 5$ are shown in Fig.~\ref{fig:result-molecules}. 

There are multiple molecules obtained from the outcomes of the MILP solver.  This was achieved using Gurobi's Solution Pool tool. Specifically, while the goal of a MILP solver is to find the globally optimal solution, other feasible solutions can be found as an algorithm progresses as an indirect output. 
We note that MILP can in this way produce a candidate pool similar to how GA can.

There are repeat molecules in the solution set.  These are solutions that have different adjacency matrices but constitute to the same molecule. In total there were 20 unique molecule-like structures found during the optimisation of the different search spaces.  The analysis of these molecules can be found in Table~\ref{table:case-study}.
\citet{Zhang23:optimizing} develop symmetry-breaking constraints, which is a route to prevent generation of repeat molecules; their work applies to MILP and also to GA approaches, and could be added to our formulations in the future.

For 12 of the 20 molecules, we found sources indicating that the molecules were experimentally observed. For the 12 observed molecules, 9 were synthesised, and their boiling points were recorded. Two of the molecules were not experimentally observed but were mentioned in papers as hypothetically possible under large pressure. Two of the molecules found during the optimisation process were also present in the training data set; all the others were not.

The absolute difference between the experimental observed boiling points and the predicted boiling points from the optimisation ranged from $9.44$ to $72.24$. The relative difference, calculated by the difference divided by the experimental boiling point, ranged from $2.7\%$--$24.2\%$. 

\begin{landscape}
\begin{table}[htp]
\centering
{\footnotesize
\begin{tabular}{cccccccccc}
\toprule
\thead{bonds} & \thead{molecule \\ length} & \thead{formula} & \thead{experimentally \\ observed} & \thead{experimental \\ $T_b$ (K)} & \thead{molecular name} & \thead{in training \\ dataset?} & \thead{objective \\ value} & \thead{predicted \\ $T_b$ (K)} & \thead{difference} \\
\midrule
{ double bonded } & 4 & C3FH5 & TRUE * & 253.15 & Allyl fluoride & FALSE & 0.02 & 313.17 & 60.02 \\
& & C2H4O2 & TRUE \cite{mardyukov20201} & & 1,1-Dihydroxyethene & FALSE & 0.78 & 361.61 & \\
& 5 & C4FH7 & FALSE & & n/a & FALSE & 0.27 & 329.79 & \\
& & C4H8O & TRUE ** & 343.15 & butanone & TRUE & 0.63 & 352.59 & 9.44 \\
& & C4H6O & FALSE & & n/a & FALSE & 0.72 & 357.90 & \\
& & C2H4O3 & TRUE ** & 375.15 & Glycolic acid & FALSE & 1.35 & 397.92 & 22.77 \\
{ triple bonded } & 4 & C3H4O & TRUE ** & 377.15 & Propargyl alcohol & FALSE & 0.57 & 348.49 & -28.66 \\
& 5 & C4H2O & TRUE \cite{araki2008laboratory} & & butadiynol & FALSE & 0.19 & 324.31 & \\
& & C4H6O & TRUE * & 329.15 & (±)-3-Butyn-2-ol & FALSE & 0.82 & 364.57 & 35.42 \\
& & C3H4O2 & FALSE & & 2-Propyne-1,1-diol & FALSE & 1.12 & 383.40 & \\
{ \thead{double \& triple \\ bonded} } & 4 & C2FH3O & TRUE *** & 283.15 & Acetyl fluoride & FALSE & 0.20 & 324.94 & 41.79 \\
& & C2F2H4 & TRUE * & 293.85 & 1,2-Difluorethane & TRUE & 0.71 & 357.30 & 63.45 \\
& & C2FH5O & TRUE ** & 366.15 & 2-Fluoroethanol & FALSE & 1.15 & 385.06 & 18.91 \\
& & CH4O3 & FALSE & & Hydroperoxymethanol & FALSE & 1.63 & 415.34 & \\
& & H2O4 & TRUE \cite{levanov2011synthesis} & & Tetraoxidane & FALSE & 1.80 & 425.95 & \\
& 5 & C4FH5 & FALSE & & 2-Fluoro-1,3-butadiene & FALSE & 0.01 & 313.17 & \\
& & C3FH5O & TRUE *** & 298.15 & 2-Fluoro-2-propen-1-ol & FALSE & 0.92 & 370.39 & 72.24 \\
& & C2H4O3 & FALSE & & Dihydroxyacetaldehyde & FALSE & 1.49 & 406.47 & \\
& & CH4O4 & FALSE \cite{bohm1997study} $^+$ & & Orthocarbonic acid & FALSE & 2.33 & 459.23 & \\
& & H2O5 & FALSE  \cite{levanov2011synthesis} $^+$ & & Pentaoxidane & FALSE & 2.39 & 463.40 & \\
\bottomrule
\end{tabular}}
\caption{Table with all experimental results and interpretation of the case study. Experimentally observed instances indicated with a star come from a database of a chemical supplier where * = Matrix Scientific, ** = Alfa Aesar, *** = Synquest. $^+ = $ hypothetical molecule in cited paper.  The predicted $T_b (K)$ is calculated by $T_b (K) = \text{mean} + \text{std} \times$ obj\_val, mean = 312.64, std = 62.98, where the mean and std are from the training data set.}
\label{table:case-study}
\end{table}
\end{landscape}

\section{Discussion}
\label{sec:discussion}

\subsection{Model Complexity}

When comparing the MILP formulation of the GCN and the MINLP formulation of the GraphSAGE network, the MINLP formulation is bi-linear, whereas the MILP formulation is linear. Linear solvers are known to be faster than non-linear solvers in general.  Another drawback for the GCN formulation is that initially, a vast amount of constraints and variable are introduced to linearise the normalisation term. This is not the case for the GraphSAGE model. 

However, the amount of constraints and variables introduced in the GraphSAGE formulation is far greater than for the GCN model for deep and wide networks. This is because for every layer $k$ the amount of constraints for the GraphSAGE model grows by $\mathcal{O}(n^2 n_k)$ whereas the GCN model only grows by $\mathcal{O}(n^2)$. For the variables we also find that the increase is of the order $\mathcal{O}(n^2 n_k)$ for GraphSAGE and $\mathcal{O}(n^2)$ for the GCN. 

While theoretical comparisons between the models has value, we also consider our empirical findings.

\subsection{Individual Experiments}
\label{subsubsection:conclusion-individual-experiments}
First we note that for very few instances the solver actually solves to optimality, for both the GCN model and the GraphSAGE model. This means that, in its current configuration, the proposed formulations can only be used to find small graphs of size $4$ in a time span of $10$ hours. Some improvements can be made to improve solving times, which will be discussed later. However, even with improvements we do not expect the search space to be able to include graphs that are multiple orders of magnitude larger than the graphs that we currently find. This means that the proposed optimisation formulations can not be used in other contexts where large graph neural networks are used, like road network modelling, or recommender systems.  This implies that that the proposed techniques, in its current formulation, should be used for small graph optimisations only, like molecule optimisation.  At the same time, the encountered computational complexity also motivates the use of inexact solvers, such as the genetic algorithm.

Second, we note that the use of a deterministic optimiser (i.e., Gurobi for MILPs) has the advantage of knowing how close one is to the actual solution, while running the algorithm, expressed by the optimality gap.  However, the experiments show that optimality gaps rapidly increase as the model becomes more complex, or as the search space includes larger graphs. When modelling some instances, this optimality gap might not be as useful anymore.  For instance, in the case that $n = 8$, the optimality gap was $26.54$ after running the $\verb|2 x 64|$ instance of the GraphSAGE model for 10 hours. The range of boiling points in the training set ranged from $145.15$--$482.05$ K. With the found objective lower bound, the optimality gap of $26.54$ implies that the solution lies in a range of approximately $675$--$2045$ K. For the set of refrigerants we could assume beforehand that the temperatures were in this boiling point range for molecules of length $8$. 

Third, we discuss the simple fact that when the number of nodes per layer increases, the solving time goes up for both the GCN and GraphSAGE. The same pattern can be seen when increasing the hidden layers per model. These results are as expected for general mixed integer linear programming formulations. As the number of layers and nodes increases, the amount of decision variables and constraints increases, making the problem more difficult to solve. The optimality gap shows similar results apart from a few exceptions as laid out in the result description section. The cases where unexpected results were seen were rerun, but resulted in similar optimality gaps, implying that the problem lay somewhere with the learned parameters of the GNN. We further explore this in the next section.

Finally, a general remark on the bounds for the GNNs. The input bound size for the MLP which comes after the pooling layer increases linearly with the amount of nodes that are in the graphs in the search space. This is because the output bounds of the feature vectors of the GNNs get summed in the pooling layer. In some cases this makes sense. Larger structures sometimes result in higher objective values, as is the case with boiling points of molecules. However, when bound are loose to start with, it amplifies this error, resulting in even larger bounds. This has a negative impact on the solving times. 

\subsection{
Comparison of the Experiments}
As discussed above, the GraphSAGE model is generally better than the GCN model in terms of solving times, and when not solved to optimality, also in terms of an optimality gap. 
There are instances where the GCN is better than the GraphSAGE model empirically.  First we note that the bounds for the GCN models are smaller than the GraphSAGE model, for equal node depth and layer width. As mentioned before, smaller bounds result in faster solving times. However, the bound difference is generally present for all instances when comparing the GCN and GraphSAGE. There thus must be another reason for these exceptions.

Training a neural network is a stochastic process. This means that training different neural networks with the same hyper-parameters does not result in the same weights and biases. We hypothesize that training different trained graph neural networks with the same configurations results in different solving times. Having different weights and biases has an impact on the bounds. In turn, we know that larger bounds have a negative impact on the solving time. We tested this hypothesis as can be seen in \ref{appendix:extra-runs}.  The same experiment for the instance $\verb|2 x 16|$ was repeated 5 times. It shows that different trained neural network parameters result in different solving times when optimised using a deterministic solver. This confirms our hypothesis. However, we find no correlation between the bounds and the solving time. We expected larger bounds to result in slower solving times, but this small test in the appendix does not confirm this hypothesis.

A final point to discuss concerning the initial experiments is the genetic algorithm baseline comparisons.  It is clear from the results that in most cases the GA is superior to the deterministic solvers in terms of finding a solution of equal quality while taking less time.  

There are three instances where the GA does not find a solution of equal quality. In these cases, the GA gets stuck in a local maximum. These instances also illustrate the main shortcoming of the GA: having no guarantee of convergence to the global optimum.

\subsection{Discussion of the Case Study}
The goal of the case study was to emulate an instance where a researcher is looking for molecules with a maximal boiling point. 
First of all, 12 of the molecules that were found were experimentally observed. Of the other 8, we able to find two which were mentioned in research as hypothetical molecules. These were able to be synthesised under very high pressure or were an unstable molecule of molecular reaction. Of the other 6, we were unable to find any mentions in literature.

We also note that only two of the 20 molecules that were found were in the original data set. We believe that this shows that a model can be trained on a particular data set and that other molecules can be found outside of that data set, of which some can be synthesised. This means there is a real-life use case for the proposed formulation in this article. Let us say a researcher wants to design a fuel with high energy storage, and low emissions. With GNNs the researcher can model these chemical properties. Using the methods proposed in this article, the researcher can make a MILP formulation of these networks, and find solutions while optimising. The researcher can use these solutions as a starting point to look for molecules with the desired properties, instead of having to start with a pool of all possible fuels.

There are two final remarks we would like to make on the found solutions. Again, these should be taken with a grain of salt, as the modelling of the chemical properties was not the main focus of this work.  However, we do see that when experimental results exist of molecules with similar input constraints and molecule length, the experimental boiling points increase as the modelled boiling points increase. This shows some validation for the modelling quality. On the other hand, we note that the mean absolute error of the trained GNN is about $6.65$.  For the found molecules, of which experimental boiling point data is available, we see that our mean absolute error is around $17.75$.  Without further exploration, we can not draw immediate conclusions from this. However, one hypothesis is that this might suggest that when modelled molecules are at the higher end of the boiling point spectrum, that the errors of the GNNs become larger.

\section{Conclusion and Outlook}
\label{section:conclusionandoutlook}\label{sec:conc}

The success of computer aided molecular design and the ubiquity of neural networks lead to the question whether one can optimally search for molecular designs, constrained by certain properties, by making use of graph neural networks.  A key barrier to using `traditional', non-graph structured networks is that they struggle to learn from non-euclidean data, whereas molecules are naturally modelled as graph-like structures, motivating the use of GNNs.

Recognising recent progress on exact formulations of non-graph neural networks as mixed integer (non-)linear programs, this work therefore formulated trained GNNs as MILP programming formulations.  These formulations can be used as surrogate models in optimisation problems.  In particular, we treated two classes of GNNs: the frequently-used GCN and the contemporary GraphSAGE.  We developed a formulation of GCN as an MINLP, and of GraphSAGE as a MILP.

In terms of accuracy, we hypothesised that the GCN would reach better model accuracy with fewer hidden layers and nodes per layer than GraphSAGE, due to the GCN model's more complex architecture.  The results (Fig.~\ref{fig:boxplots}), do not support this hypothesis.  For four of the seven configurations (hidden layers $\times$ nodes) the GraphSAGE model has a better median validation MSE while running for 20 iterations, and for four of the seven configurations, the GraphSAGE model has a lower minimum validation error than the GCN model.  Overall, with the hyper-parameters tested, we achieved similar model accuracy for both the GCN and GraphSAGE model, even with the same number of hidden layers and nodes per layer. 

In terms of solving speed, we hypothesised was that the GraphSAGE MILP formulations would be faster than the GCN MINLP formulations because linear solvers are faster than non-linear solvers (in this case, bi-linear).
The results (Tables~\ref{tab:results-GCN} and~\ref{tab:results-GraphSAGE})
find that the GraphSAGE model was generally faster (four of the six solved instances).  The optimality gaps also seem to suggest that if the experiments were ran for longer the GraphSAGE would generally solve to optimality first.  This is because for all but three configurations (12 out of 15 early terminated cases) the GraphSAGE model had a smaller optimality gap than the GCN model. 
Overall, there is evidence to suggest the MILP formulation of the GraphSAGE model solves to optimality faster than the MINLP formulation of the GCN model, with similar model accuracy. This is because our trained model accuracy is about the same and sometimes better for the GraphSAGE model compared to the GCN model, for models with similar hidden layers and number of nodes, combined with the fact that the GraphSAGE model often solves to optimality faster with similar configurations.

Our final contribution was to apply the MI(N)LP formulations to a case study of optimising the boiling points of molecules.  The case study successfully derived a set of optimal molecules, given constraints on the design space.  Of the 20 molecules derived, 12 were found were experimentally observed.  Of the other eight, the literature notes two as hypothetical molecules.  These were able to be synthesised under very high pressure or were an unstable molecule of molecular reaction.  The remaining six molecules appear to be novel; their chemical feasibility in practice would be tested in vitro studies.

Our work opens up several prominent research directions.  First, the models themselves have potential for improvement with stronger bound tightening techniques, and we think techniques for tightening MLP MILP formulations can be applied to GNN MILPs also.  Going beyond feasibility-based bound tightening, optimisation-based bound tightening techniques (OBBT), and the combined technique of \citet{wang2021acceleration}. Second, using GNNs in CAMD, there is opportunity in increasing the training set size and using more of the learnable features, and reconsidering linearise structures and the input constraints.  We underline that \emph{training} of the GNNs was not the main focus of the current work, and that we applied our novel formulations to GNN as a case study.
Third, when deterministic optimisation is not the main priority for researchers, our straightforward GA for optimising trained GNNs already shows promise.

\subsubsection*{Funding}
This work was partially supported EU Horizon 2020, grant number 952215 (TAILOR).

\iffalse
\appendix
\include{Appendices/Appendix-ReLU-NN-as-MILP}

\include{Appendices/A1-MINLP-of-GNN-Appendix}

\include{Appendices/Appendix-Molecule-Constraints}

\include{Appendices/Appendix-GA}

\include{Appendices/Appendix-Initial-Experiments}

\include{Appendices/A2-Extra-Results-Appendix}
\fi

\bibliographystyle{elsarticle-num-names}
\bibliography{references,more-refs}

\end{document}